\documentclass[12pt]{amsart}
\usepackage{amssymb,amsmath}
\usepackage{andy}
\usepackage{amsthm}
\usepackage{pdfsync}

\newtheorem{thm}{Theorem}[section]
\newtheorem{prop}[thm]{Proposition}
\newtheorem{lem}[thm]{Lemma}
\newtheorem{cor}[thm]{Corollary}

\theoremstyle{definition}

\theoremstyle{remark}
\newtheorem{rem}[thm]{Remark}

\newcommand{\Top}{\mathcal T_{\ell, \aleph}}
\newcommand{\AdjTop}{\mathcal T^*_{\ell, \aleph}}

\begin{document}

\title[ Div-Curl for Higher Order] {On Div-Curl for Higher Order}
\author[Lanzani and Raich]{Loredana Lanzani$^*$ and Andrew S. Raich$^{**}$}
\thanks{To appear in {\em Advances in Analysis: the Legacy of E. M. Stein}, Princeton U. Press.}
\thanks{$^*$ Supported by a National Science Foundation IRD plan, and in part by  awards DMS-0700815 and DMS-1001304}
\thanks{$^{**}$ Supported in part by the National Science Foundation, award DMS-0855822}
\thanks{THIS IS AN ELECTRONIC REPRINT OF THE ORIGINAL ARTICLE THAT APPEARED IN {\em ADVANCES IN ANALYSIS: THE LEGACY OF E. M. STEIN}, PRINCETON U. PRESS (2013) ISBN: 9780691159416. THIS REPRINT DIFFERS FROM THE ORIGINAL IN PAGINATION AND TYPOGRAPHICAL DETAIL}

\address{Dept. of Mathematics \\
University of Arkansas Fayetteville, AR 72701}

\email{lanzani@uark.edu,\ araich@uark.edu}
\maketitle

\centerline{\em For Eli}

\begin{abstract}
 We present new examples of complexes of differential operators of order $k$ (any given positive integer) that satisfy
div-curl and/or $L^1$-duality estimates. 
\end{abstract}

\section{Introduction}\label{S:introduction}

In 2004 Stein and the first named author discovered a connection \cite{LS} between the celebrated Gagliardo-Nirenberg inequality \cite{G}-\cite{N} for functions
\begin{equation}\label{E:GN-classical}
\|f\|_{L^r(\mathbb R^n)} \leq C \|\nabla f\|_{L^1(\R^n)},\quad r=n/(n-1)
\end{equation}
and a recent estimate
   of Bourgain and Brezis \cite{BB2} for divergence-free vector  fields 
  as proved by Van Schaftingen \cite{VS1}
 \begin{equation}\label{E:BB-first}
 \|Z\|_{L^r(\mathbb R^n))}\leq C\|\mbox{Curl}\, Z\|_{L^1(\mathbb R^n)},\quad r=n/(n-1),\quad  \mbox{div}\,Z=0
 \end{equation}
Such connection is 
provided by
  the exterior derivative operator
   acting on differential forms on $\mathbb R^n$
with (say) smooth and compactly supported coefficients
\begin{equation}\label{E:ext-der}
d:\  \Lambda_q(\mathbb R^n) \to \Lambda_{q+1}(\mathbb R^n),\quad
 0\leq q\leq n
\end{equation}
 It was proved in \cite{LS} that the inequality
\begin{equation}\label{E:GN}
\|u\|_{L^r(\R^n)}\leq C (\|du\|_{L^1(\R^n)}+\|d^*u\|_{L^1(\R^n)}),\quad r=n/(n-1)
\end{equation} 
holds for any form $u$ of degree $q$ other than $q=1$ (unless $d^*u=0$) and $q= n-1$ (unless $du=0$). Note that 
\eqref{E:GN-classical} is the case $q=0$, whereas \eqref{E:BB-first} is the case $q=1$ specialized to $d^*u=0$.
\medskip

\noindent Since those earlier results 
 div/curl-type 
 phenomena have been 
 studied
   both in the Euclidean and non-Euclidean settings  
\cite{Am}, \cite{BV}, \cite{HP1}, \cite{HP2}, \cite{M}, \cite{MM}, \cite{Mi}
  \cite{VS4},
 \cite{CV}, \cite{Y}.
 In \cite{VS2} and the recent works \cite{BB3}, \cite{VS3}, \cite{VS5},
  differential conditions of higher order have been considered for the first time in such context.
   (By contrast, the exterior derivative in \eqref{E:ext-der} 
  is defined in terms of differential conditions of order $1$.)
 \medskip
 
   The goal of the present paper is to produce a new class of differential operators of order $k$ 
   (where $k$ is any given positive integer) that satisfy an appropriate analogue of
   \eqref{E:GN} and contain the operators introduced in \cite{BB3}, \cite{VS2} and \cite{VS3}; since the conditions
      \begin{equation}\label{E:complex}
   d\circ d=0;\qquad d^*\circ d^*=0
   \end{equation}
   play an important role in the proof of \eqref{E:GN}, the new operators
   should satisfy \eqref{E:complex} as well.
   We achieve this goal in a number of 
     ways, beginning with:
   \begin{thm}\label{T:LR-BB} If $u\in C^\infty_q(\mathbb R^n)$ has compact support, then
 \begin{equation}\label{E:GN-k}
\|u\|_{W^{k-1, r}}\leq C (\|\mathcal T u\|_{L^1}+\|\mathcal T^*u\|_{L^1}),\quad r=n/(n-1)
\end{equation}
whenever $q$ is neither $1$ (unless $\mathcal T^*u=0$)  nor $n-1$ (unless $\mathcal T u=0$), where
\begin{equation}\label{E:LR-BB-op}
\mathcal T u := 
\sum\limits_{|L|=q+1}
\!\!\bigg(
  \sum\limits_{\atopp{
  |I|=q
  }{j=1,\ldots, n}}
  \!\!\!
  \epsilon_L^{j I}\
  \frac{\dee^k u_I}{\dee x_j^{k}}
  \bigg)dx^L
\end{equation}
 \end{thm}
 Here and in the sequel, $W^{a, p}(\mathbb R^n)$ denotes the Sobolev space consisting of
 $a$-times differentiable functions in the Lebesgue class $L^p(\mathbb R^n)$ (and 
 $W^{a, p}_q(\mathbb R^n)$ will denote the space of $q$-forms with coefficients in
  $W^{a, p}(\mathbb R^n)$), while  $\epsilon^{AB}_C\in\{-1, 0, +1\}$ is the sign
  of the permutation that carries the ordered set
   $AB=\{a_1,\ldots, a_\ell, b_1,\ldots,b_q\}$ to the label $C=(c_1,\ldots,c_{\ell+q})$, if these have identical content, and is otherwise zero. 
  Note that when $k=1$ then $\mathcal T = d$ and
 inequality \eqref{E:GN-k} is indeed \eqref{E:GN}.

 \medskip
  
Another such complex,
again involving a differential condition of order $k\geq 1$,
is obtained
by embedding $\mathbb R^n$ isometrically in a larger space $\mathbb R^N$.
(The choice of ``inflated'' dimension $N$ will be discussed later.)
The resulting operators act on ``hybrid $\mathbb R^n$-to-$\mathbb R^N$'' spaces of forms
 whose coefficients
 are trivial extensions to $\mathbb R^N$ of functions defined on
$\mathbb R^n$; to distinguish such spaces from the classical Sobolev spaces 
$W^{a, p}_q(\mathbb R^N)$
 (to which they are by necessity transversal) we will use the notation 
$$\widetilde W^{a, p}_q(\mathbb R^N),\qquad 0\leq q\leq N,$$
and we will 
write
$$
\widetilde C^{\infty, c}_q(\mathbb R^N),\qquad 0\leq q\leq N
$$
to indicate a dense subspace of smooth ``compactly supported'' forms.
These operators, which we denote $\widetilde T_{1, \aleph}$, map
$$
\widetilde T_{1, \aleph}: \widetilde C^{\infty, c}_q(\mathbb R^N) \to
\widetilde C^{\infty, c}_{q+1}(\mathbb R^N),\qquad 0\leq q\leq N
$$
The label $\aleph$ refers to a choice of an ordering for the set of all $k$-th order derivatives in $\mathbb R^n$, and so in practice we define a finite family 
$\{\widetilde T_{1, \aleph}\}_{\aleph}$ of such complexes. 
(We use the subscript ``$1$''    in $\widetilde T_{1, \aleph}$ to specify that 
$\widetilde T_{1, \aleph}$ maps $q$-forms
to $(q+1)$-forms, a distinction that will be needed later on).
The explicit definition of $\widetilde T_{1, \aleph}$ will be given in the next section; what matters here 
is that
these operators satisfy a more general version of \eqref{E:GN-k}
 in the sense that 
the following inequality
 implies 
  \eqref{E:GN-k}
 but the converse is not true:
\begin{thm}\label{T:LR-VS-1}
If $U\in \widetilde C^\infty_q(\mathbb R^N)$ has compact support, then
 \begin{equation}\label{E:LR-BB-N}
\|U\|_{\widetilde W^{k-1, r}}\leq C (\|\widetilde T_{1, \aleph} U\|_{\widetilde L^1}+\|\widetilde T^*_{1, \aleph}U\|_{\widetilde L^1}),\quad r=n/(n-1)
\end{equation}
whenever $q$ is neither $1$ (unless $\widetilde T^*_{1, \aleph} U=0$) nor $N-1$ 
(unless $\widetilde T_{1, \aleph} U=0$).
\end{thm}
Theorem \ref{T:LR-BB}  recaptures
 an $L^1$-duality estimate 
of Bourgain and Brezis:
\begin{thm}[\cite{BB3}]\label{T:BB-thm-3} 
 Let $k\ge 1$. For every vector field 
$\displaystyle{u\in L^1(\mathbb R^n;\ \mathbb R^n)}$ 
 if
\begin{equation}\label{E:BB-div-k}
\sum\limits_{j=1}^n\frac{\dee^k u_j}{\dee x_j^k} =0
\end{equation}
in the sense of distributions, then
\begin{equation}\label{E:BB-ineq-div-k}
\bigg|\,\int\limits_{\mathbb R^n}u_j\cdot h_j\,\bigg|\leq C\|u\|_{L^1}\|\nabla h\|_{L^n},\qquad j=1,\ldots, n
\end{equation}
for any $\displaystyle{h\in (W^{1, n}\cap L^\infty)(\mathbb R^n;\ \mathbb R^n)}$,
where the constant $C$ only depends on the dimension of the space
 $n$ and on the order $k$.
 \end{thm}
On the other hand, Theorem \ref{T:LR-VS-1}  was motivated by
a recent result of van Schaftingen:
 \begin{thm}[\cite{VS3}]\label{T:VS-thm-4} 
 Given $k\ge 1$ and $n\ge 2$, let 
  \begin{equation}\label{E:M-def}
m:= 
\binom{n-1 +k}{k}
\end{equation}
For any vector field 
$\displaystyle{g =(g_\alpha)_{\alpha\in\mathcal S (n, k)}\in L^1(\mathbb R^n;\ \mathbb R^m)}$ 
 if
\begin{equation}\label{E:VS-div-k}
\sum\limits_{\alpha\in\mathcal S(n, k)}\frac{\dee^k g_{\alpha}}{\dee x^\alpha} =0
\end{equation}
in the sense of distributions, then
\begin{equation}\label{E:VS-ineq-div-k}
\bigg|\,\int\limits_{\mathbb R^n}g_\alpha\cdot h_\alpha\,\bigg|\leq C\|g\|_{L^1}\|\nabla h\|_{L^n},\qquad \alpha\in\mathcal S(n, k)
\end{equation}
for any $\displaystyle{h
\in (W^{1, n}\cap L^\infty)(\mathbb R^n;\ \mathbb R^m)}$,
where the constant $C$ only depends on the dimension of the space
 $n$ and on the order $k$.
 \end{thm}
 Here 
 $\mathcal S (n, k)$ denotes the set of $k$-multi-indices in $\mathbb R^n$:
  \begin{equation}\label{E:S(n,k)}
    \mathcal S (n, k) =
 \bigg\{ \alpha = (\alpha_1,\ldots,\alpha_n)\,\bigg| \,
 \alpha_t\in \{0, 1,\ldots, k\},\,
 \sum\limits_{t=1}^n\alpha_t
 = k\bigg\}.
 \end{equation}

A key ingredient in the proof of Theorems \ref{T:LR-BB} and \ref{T:LR-VS-1} is the fact that the Hodge laplacians for these operators, namely
$$
\square_{\mathcal T}:= \mathcal T\mathcal T^* + \mathcal T^*\mathcal T:\ \ 
C^{\infty, c}_q(\mathbb R^n)\to C^{\infty, c}_q(\mathbb R^n),\ \ 0\leq q\leq n
$$
and
$$
\widetilde \square_{1, \aleph}:= \widetilde T_{1, \aleph}\widetilde T^*_{1, \aleph} +
\widetilde T_{1, \aleph}^*\widetilde T_{1, \aleph}:\ \  \widetilde  C^{\infty, c}_q(\mathbb R^N)\to \widetilde  C^{\infty, c}_q(\mathbb R^N),\ \ 0\leq q\leq N
$$
satisfy a uniform Legendre-Hadamard condition
 which in turn yields elliptic estimates. 
\medskip

Rather surprisingly, it turns out that 
in fact there is 
a larger class
of such operators, mapping
$$
\widetilde T_{\ell, \aleph}: \widetilde C^{\infty, c}_{q} (\mathbb R^N) \to
\widetilde C^{\infty, c}_{q + \ell} (\mathbb R^N),\quad 0\leq q\leq N
$$
where the label $\ell$ now runs over all the elements of what we call   
{\em the set of admissible degree increments}, which is a subset of 
$\{1,\ldots , k\}$ determined by $n$ (the dimension of the source space) and $k$ (the order of differentiation): for any $n\geq 2$ and $k\geq 2$, the set of admissible degree increments
 contains at least two distinct elements: $\ell =1$ (discussed earlier) and also
  $\ell =k$. Each admissible degree increment in turn determines an
  ``inflated dimension'' $N$
  (in particular
  $N$ will change with $\ell$).
However the situation for $\ell \neq 1$ differs from the case $\ell=1$ in two important respects: the crucial condition \eqref{E:complex} will hold only for 
odd $\ell$, and if $\ell \neq 1$ the Hodge laplacian for 
$\widetilde T_{\ell, \aleph}$ will fail to be uniformly elliptic (even for $\ell$ odd): as a result
there is no analog of \eqref{E:GN-k}.
Instead, we show that for any admissible degree increment
(thus also for $\ell =1$), the operators 
$\widetilde T_{\ell, \aleph}$ satisfy $L^1$-duality estimates that are similar in spirit, and indeed are equivalent to \eqref{E:VS-ineq-div-k};
 see Theorem \ref{P:LR-inflated-dim} for the precise statement.
\medskip

A further class of operators which contains our very first example 
$\mathcal T$, see \eqref{E:LR-BB-op},
can be defined 
 in terms of 
$\widetilde T_{\ell, \aleph}$ and of the 
aforementioned
 embedding: 
$\mathbb R^n\to\mathbb R^N$.
Such operators map
$$
\mathcal T_{\ell, \aleph}: C^{\infty, c}_q(\mathbb R^n)\to 
C^{\infty, c}_{q+\ell}(\mathbb R^n)
$$
and satisfy div-curl 
and/or 
$L^1$-duality estimates
that are stated solely in terms
 of the 
  source space $\mathbb R^n$
 rather than the ``hybrid $\mathbb R^n$-to-$\mathbb R^N$'' spaces 
$\widetilde L^p_q(\mathbb R^N)$ and 
$\widetilde W^{a, p}_q(\mathbb R^N)$,
see Theorem \ref{P:LR-source-smallest} and \eqref{E:f-q}.

(Of course, if $\ell\neq 1$ such operators are non-trivial 
only for $n\geq \ell$.)
\medskip

We need to explain the reason for our choice to keep track, through the label $\aleph$, of the orderings of $\mathcal S(n, k)$: this has to do with the notion of invariance.
 One would like to know whether the identity
\begin{equation}\label{E:invariance-big}
\widetilde T_{\ell, \aleph}\,\Psi^* F = \Psi^*\,\widetilde T_{\ell, \aleph} F
\end{equation}
holds for any $F\in \widetilde C^{\infty}_q(\mathbb R^N)$ and for some non-trivial class of diffeomorphisms 
$$
\Psi: \mathbb R^N\to \mathbb R^N
$$
of class $C^{k+1}$:
it is in this context that the choice of $\aleph$ may be relevant. In the case $k=1$ 
 our construction
gives $N=n$ with $\aleph$
spanning  the set of all permutations of $\{1, \ldots, n\}$, and since $k$ is 1 there is only the admissible degree increment $\ell =1$. As a result, for $k=1$ we have
$$
\widetilde T_{1, \aleph}U =
\sum\limits_{|L|=q+1}
\!\!\bigg(
  \sum\limits_{\atopp{
  |I|=q
    }{j=1,\ldots, n}}
  \!\!\!
  \epsilon_L^{\aleph(j) I}\
  \frac{\dee U_I}{\dee x_j}
   \bigg)dx^L,\quad
 \aleph\in \Sigma (1,\ldots, n).
$$

In particular one has
$$
\widetilde T_{1, \aleph_0} = \mathcal T = d
$$
for exactly one permutation $\aleph_0$ (the identity) which therefore determines 
an invariant operator. On the other hand it is easy to check that for any $\aleph\neq\aleph_0$
the operators $\widetilde T_{1, \aleph}$ fail to be invariant.
   \smallskip
   
No such phenomenon exists for $k\geq 2$: there is no choice of $\aleph$ 
(nor $\ell$) that makes $\widetilde T_{\ell, \aleph}$ invariant and  \eqref{E:invariance-big} fails even in the case
when $\Psi$ originates from a rotation of $\mathbb R^n$. 
It can be verified that $\mathcal T_{\ell, \aleph}$, too, is not invariant because if $k\geq 2$ the identity
\begin{equation}\label{E:invariance-small}
\mathcal T_{\ell, \aleph}\,\psi^* = \psi ^*\,\mathcal T_{\ell, \aleph} 
\end{equation}
  fails  for any $\ell$ and for any $\aleph$, already for $\psi$ a rotation of $\mathbb R^n$.   
 \medskip
  
  Finally, we point out that our results can be rephrased in terms of the canceling and cocanceling conditions 
 of \cite{VS4}: within that framework,
   our results provide new classes of differential operators of 
 {\em arbitrary} order that are canceling and/or cocanceling, with the size of the admissible degree increments acting as an indicator of the canceling property. See the remarks in Section \ref{S:open-problems}.

\medskip

This paper is organized as follows: in Section \ref{S:statements}
we introduce the notion of admissible degree increment,
  we describe the ``hybrid $\mathbb R^n$-to-$\mathbb R^N$''
  Sobolev spaces $\widetilde W^{a, p}_q(\mathbb R^N)$ in term of the embedding, and we define the operators
$\widetilde T_{\ell, \aleph}$ and $\mathcal T_{\ell, \aleph}$ and discuss their basic properties (adjoints; uniform ellipticity). The $L^1$-duality estimates  for $\widetilde T_{\ell, \aleph}$ and for
$\mathcal T_{\ell, \aleph}$ are stated in Theorems 
\ref{P:LR-inflated-dim} and \ref{P:LR-source-smallest}, and  the precise 
 statements of \eqref{E:LR-BB-N} and of
 \eqref{E:GN-k}
 are given  in Theorem \ref{T:LR-thm-Z} and in \eqref{E:f-q}.
 All the proofs are deferred to Section \ref{S:Proofs}. Section
  \ref{S:open-problems} contains some remarks and a few questions.

\medskip

\subsection{Notation}\label{Ss:prelim}
As customary, we let $\Lambda_q(\mathbb R^n)$ denote the space of $q$-forms:
   \begin{equation}\label{E:lambda-q}
 \Lambda_q(\mathbb R^n) =
 \bigg\{f= \!\!\!\sum\limits_{I\in\mathcal I(n, q)}\!\! f_I\,dx^I\ \big|\ f_I:\mathbb R^n\to \mathbb R\bigg\},\quad 0\leq q\leq n
 \end{equation}
  where $\mathcal I(n, q) $ denotes the set of $q$-labels for $\mathbb R^n$:
 \begin{equation}\label{E:I(n,q)}
 \mathcal I(n, q) =
 \big\{ I = (i_1,\ldots,i_q)\,|\, i_t\in \{1,\ldots, n\}, \ i_t<i_{t+1}\big\}
 \end{equation}
 and 
 $$dx^I = dx_{i_1}\wedge\cdots\wedge dx_{i_q}.$$
 
 When $q=n$ the expression above is the volume form and we use the notation $dV$. We will regard the label set $\mathcal I(n, q)$ as canonically ordered (alphabetical ordering).
  Letting
 $$
 i:\mathbb R^n \to \mathbb R^N
 $$ 
 denote the 
   isometric embedding mentioned above and defined in \eqref{E:canonical-embed}, 
  the ``hybrid $\mathbb R^n$-to-$\mathbb R^N$'' subspace of $\Lambda_q(\mathbb R^N)$ (consisting of those
 $q$-forms whose coefficients are trivial extensions to $\mathbb R^N$ of functions
 defined on $\mathbb R^n$)
  is more
 precisely described as follows 
  \begin{equation}\label{E:def-tilde-q-N}
   \widetilde\Lambda_q(\mathbb R^N):=
   \bigg\{F=\sum\limits_{I\in\mathcal I(N, q)}\!\!\!\!F_I\,dz^I
    \ \bigg|\ F_I\,\circ\, i\,\circ\,\pi = F_I\bigg\},
   \ 0\leq q\leq N
   \end{equation}
   where $$\pi: \mathbb R^N\to \mathbb R^n$$ is chosen so that  
\begin{equation}\label{E:pi-i}
(\pi\circ i) (x)\ =\ x,\quad \text{for all}\quad x\in\mathbb R^n
\end{equation}
 As a result  the reverse composition
\begin{equation}\label{E:pi-i-pr}
i\circ \pi: \mathbb R^N\to \mathbb R^N
\end{equation} 
is a projection.  

We will denote the Hodge-star operators 
for each of $\Lambda_q(\mathbb R^n)$ and $\Lambda_q(\mathbb R^N)$ respectively by 
$*_n$ and $*_N$; note that we have
\begin{equation}\label{E:star-tilde}
*_N:\widetilde\Lambda_q(\mathbb R^N)\to \widetilde\Lambda_{N-q}
(\mathbb R^N),\quad 0\leq q\leq N.
\end{equation}

\bigskip

\section{Statements}\label{S:statements}
\medskip

\subsection{Admissible degree increments}\label{Ss:adm-deg}
 Given three integers: 
 \medskip
 
 \begin{itemize}
\item[i.\quad ] $n\ge 2$ (the dimension of the source space), 
\medskip

\item[ii.\quad] $k\ge 1$ (the order of the differential condition), and 
\medskip

\item[iii.\quad] $1\leq \ell\leq k$,
 \end{itemize}
 \medskip
 
\noindent  we say that $\ell$ is an {\it admissible degree increment for the pair
   $(n, k)$} if and only if the
  polynomial equation
  \begin{equation}\label{E:embedding-dim}
  \binom{N}{\ell} = \binom{n-1+k}{k}
\end{equation}
has a solution $N$ that satisfies the following two conditions:
\begin{equation}\label{E:constraint-dim}
N\in\mathbb Z^+,\ \ N \ge n-1+\ell\, .
\end{equation}
Note that the pair $(n, 1)$ (that is, $k=1$) has {\it exactly one} admissible
degree increment, namely $\ell =1$, and in this case equation 
\eqref{E:embedding-dim} has the unique solution: $N=n$.
On the other hand, for $k\ge 2$ any pair $(n, k)$ will 
have {\it at least two} admissible degree increments ($\ell =1, k$) and possibly more, for instance: the pair $(n, k)=(2, 9)$ has (exactly) four admissible degree increments, namely $\ell =1, 2, 3, 9$; similarly, the
pair $(n, k):=(2, 29)$ has  (at least)
$\ell = 1, 2, 29$. For any admissible degree increment, we consider the embedding
  \begin{equation}
  \label{E:canonical-embed-1}
   i: \mathbb R^n\to \mathbb R^N
   \end{equation}
   defined as follows
   \begin{equation}\label{E:canonical-embed}
    i(x_1,\ldots, x_n)
   = (z_1,\ldots, z_N) := (x_1,\ldots, x_n, 0,\ldots, 0)
   \end{equation}
    where $N = N(n, k, \ell)$ is as in \eqref{E:embedding-dim} and
\eqref{E:constraint-dim}.
We let $i$ also 
denote  the embedding
of $k$-multi-indices 
        \begin{equation*}
   i: \mathcal S(n, k) \to 
     \mathcal S (N, k)
   \end{equation*}
   that is canonically induced by \eqref{E:canonical-embed}, namely
   \begin{equation}\label{E:multi-index-embed-1}
    i (\alpha_1,\ldots, \alpha_n):=
         (\alpha_1,\ldots,\alpha_n,0,\ldots,0)\in \mathcal S(N, k)
      \end{equation}
      and adopt the notation
      \begin{equation}\label{E:multi-index-embed-2}
      i\mathcal S(n, k) := \{i\alpha\, |\ \alpha\in \mathcal S(n, k)\}
      \subsetneq
      \mathcal S(N, k)
      \end{equation}
  We have
\begin{equation}\label{E:constraint}
\left|i\mathcal S(n, k)\right| = 
m
\end{equation}
with $m = m(n, k)$ as in \eqref{E:M-def}, and so there are $m!$-many
distinct orderings of $i\mathcal S(n, k)$. By the definition of $N$ 
the set
of labels
 $\mathcal I(N, \ell)$  
  also has
cardinality $m$
  and 
   we will think of each ordering of $i\mathcal S(n, k)$ as a
one-to-one correspondence
\begin{equation}\label{E:correspond}
\aleph:\  i\mathcal S(n, k)\ \to 
 \  \mathcal I (N, \ell);\qquad
 \aleph^{-1}:\,  \mathcal I (N, \ell)\ \to\ i\mathcal S(n, k)
\end{equation}
\bigskip

\subsection{Hybrid Function spaces}\label{Ss:function}
Given an  integer $a\geq 0$ 
and given $p, p'\geq 1$ such that $1/p+ 1/p'=1$, we first set ($q=0$)

\begin{equation}\label{E:Lebesgue-big}
\widetilde L^p(\mathbb R^N):=\{F:\mathbb R^N\to \mathbb R\ |\ F\circ i\circ \pi = F,\ \ F\circ i\in L^p(\mathbb R^n, dV)\}
\end{equation}
where $i$ is as in \eqref{E:canonical-embed} and 
\begin{equation}\label{E:canonical-proj}
\pi (z_1,\ldots, z_n,\ldots, z_N)=(x_,\ldots, x_n):= (z_1,\ldots, z_n)
\end{equation}
satisfies \eqref{E:pi-i},
and
\begin{equation*}
\widetilde W^{a, p}(\mathbb R^N)\! :=\!
\Big\{F\!:\!\mathbb R^N\!\to\mathbb R\, |\, \dee^\lambda F\in 
\widetilde L^{p}(\mathbb R^N),\,\lambda\in\mathcal S(N, s),\, 0\leq s\leq a\Big\}.
\end{equation*}
Note that if $F\in \widetilde W^{a, p}(\mathbb R^N)$ then
it follows from
\begin{equation*}\label{E:sobolev-deriv-0}
F\circ i\circ \pi = F
\end{equation*}
that 
\begin{equation}\label{E:sobolev-deriv-1}
 \frac{\dee F}{\dee z_t} =0, \quad \mbox{for any}\ t=n+1, \ldots, N, \ 
\end{equation}
which in turn grants
\begin{equation}\label{E:sobolev-deriv-2}
 \frac{\dee^s F}{\dee z^\lambda} =0
 \end{equation}
 for any $ 1\leq s\leq a$
 and for any $\lambda\in\mathcal S(N, s)\setminus i\mathcal S(n, s)$, 
 so that these spaces are more precisely described as follows
 \begin{equation*}
\widetilde W^{a, p}(\mathbb R^N)\! =\!
\Big\{F\!:\!\mathbb R^N\!\to\mathbb R\, |\, \dee^{i\beta} F\in 
\widetilde L^{p}(\mathbb R^N),\,\beta\in\mathcal S(n, s),\, 0\leq s\leq a\Big\}.
\end{equation*}

As customary, these definitions are extended to forms
$F\in \widetilde L^p_q(\mathbb R^N)$ 
of any degree $0\leq q\leq N$ 
$(\mbox{resp.} \ F\in \widetilde W^{a, p}_q(\mathbb R^N), 0\leq q\leq N)$
by requiring that
$$
F = \!\!\sum\limits_{I\in\mathcal I(N, q)}\!\!\!\!F_I\,dz^I
\quad \mbox{has }
\ F_I\in \widetilde L^p(\mathbb R^N) \ \
\mbox{(resp.}\ F_I\in \widetilde W^{a, p}(\mathbb R^N)\mbox{)}\quad 
$$
$\mbox{for any}\ I\in\mathcal I(N, q)$.
We observe for future reference that identity \eqref{E:pi-i} grants
\begin{equation}\label{E:sobolev-deriv-3}
 \left(\frac{\dee^s F}{\dee z^{i\beta}}\right)\circ i = 
 \frac{\dee^s (F\circ i)}{\dee x^\beta}
 \end{equation}
 for any 
 $\beta\in \mathcal S(n, s)$ and for any 
 $F\in \widetilde W^{a, p}(\mathbb R^N)$.
\bigskip

\begin{lem}\label{L:hilbert}
For any $0\leq q\leq N$; for any $p\geq 1$ and for any  integer $a\geq 1$ the following properties hold:
\medskip

\begin{itemize}

\item[{\em {i}.}\quad ] $\widetilde L^2_q(\mathbb R^N)$ is a Hilbert space with inner product
\begin{equation}\label{E:pairing-big}
\langle F, G\rangle_{\widetilde L}\ := 
\int\limits_{\mathbb R^n} *_n\,i^*\!*_N(F\wedge *_N G)
\end{equation}
\medskip

\item[{\em ii.}\quad ]$\widetilde L^p_q(\mathbb R^N)$ is a Banach space with norm
\begin{equation}\label{E:Leb-norm-big}
\|F\|_{\widetilde L^p_q}:=
\left(\
\int\limits_{\mathbb R^n} *_n\,i^*\big(\!*_N(F\wedge *_N F)\big)^{\!p/2}
\right)^{\!\!\!1/p}
\end{equation}
\medskip

\item[{\em{iii.}}\quad ] $\widetilde W^{a, 2}_q(\mathbb R^N)$ is a Hilbert space
with inner product
\begin{equation}\label{E:inner-sob-big}
(\,F, G\,)_{\widetilde W}\ := \
\sum\limits_
{\atopp{0\leq s\leq a}{\beta \in\mathcal S(n, s)}}
\langle D^{i\beta}F, D^{i\beta}G\rangle_{\widetilde L}\, 
\end{equation}
where we have set
\begin{equation}\label{E:def-deriv-form}
D^{i\beta}F:= \sum\limits_{I\in\mathcal I(N, q)}\!\!\!
\left(\dee^{i\beta}F_I\right) dz^I
\end{equation}
\medskip

\item[{\em iv.}\quad ] $\widetilde W^{a, p}_q(\mathbb R^N)$ is a Banach space with norm
\begin{equation}\label{E:Sob-norm-big}
\|H\|_{\widetilde W^{a, p}_q}:=
\Bigg(
\sum\limits_
{\atopp{I\in\mathcal I(N, q)}
{\atopp{\beta\in\mathcal S(n, s)}{0\leq s\leq a}}}
  \!\!\!
  \|\dee^{i\beta} H_I\|^p_{\widetilde L^p}
\Bigg)^{\!1/p}
\end{equation}
\medskip

\item[{\em{v.}}\quad ]The set
\begin{equation}\label{E:dense-big}
\widetilde C^{\infty, c}_q(\mathbb R^N) \ :=\ 
\left\{
F=\!\!\sum\limits_{I\in\mathcal I(N, q)}\!\!\!\!F_I\,dz^I
\in \widetilde\Lambda_q(\mathbb R^N)\, \Bigg|\, F_I\circ i\in 
C^\infty_c(\mathbb R^n)
\right\}
\end{equation}
is dense in $\widetilde L^p_q(\mathbb R^N)$  
(resp. $\widetilde W^{a, p}_q(\mathbb R^N)$) 
with respect to the norm
\eqref{E:Leb-norm-big}  (resp. \eqref{E:Sob-norm-big}).
\end{itemize}
\end{lem}

A sequence $\{\Phi_j\}_j\subset \widetilde C^{\infty, c}_q(\mathbb R^N)$ is said to converge in the sense of the space 
$\widetilde{\mathcal D}_q (\mathbb R^N)$ to 
$\Phi\in \widetilde C^{\infty, c}_q(\mathbb R^N)$, see \cite{A}, provided the following conditions are satisfied:
\begin{itemize}
\item[{\em i.}] There is a set $K\Subset \mathbb R^n$ such that 
$$\mbox{Supp(}(\Phi_j-\Phi)\circ i\mbox{)}\subset K \ \mbox{for each}\ j$$
\item[{\em ii.}] For any $0\leq s<\infty$ and for each $\beta\in\mathcal S(n, s)$ we have
$$\displaystyle{\lim\limits_{j\to\infty}
\frac{\dee^\beta ((\Phi_{j})_I\circ i)}{\dee x^\beta} =
\frac{\dee^\beta (\Phi_{I}\circ i)}{\dee x^\beta}\quad 
\mbox{uniformly in}\ K
}$$
 for any $I\in\mathcal I(N, q)$. 
\end{itemize}
There exists a locally convex topology on the vector space
 $\widetilde C^{\infty, c}_q(\mathbb R^N)$ with respect to which a linear functional $\mathcal L$
is continuous if, and only if, $\mathcal L (\Phi_j)\to \mathcal L (\Phi)$ whenever
$\Phi_j\to\Phi$ in the sense of the space 
$\widetilde{\mathcal D}_q(\mathbb R^N)$.
\smallskip

For any $1\leq p, p'<\infty$ with $1/p + 1/p' =1$, the dual space
 $\widetilde W^{-a, p'}_q(\mathbb R^N)$ of 
 $\widetilde W^{a, p}_q(\mathbb R^N)$ is identified 
 (in the usual fashion, see e.g., \cite[III.3.8 -- III.3.12]{A}) with a 
 closed subspace of the Cartesian product
 $$
 (\widetilde L^{p'}(\mathbb R^N))^{\binom{N}{q} M}
 \quad
 \mbox{where}\quad 
 M=M(n, a):=\sum\limits_{j=0}^a\binom{n-1+j}{j}
 $$
and from this it follows that for any $F\in\widetilde W^{a, p}_q(\mathbb R^N)$ and 
$G\in \widetilde W^{-a, p'}_q(\mathbb R^N)$ we have 
\begin{equation}\label{E:holder-big}
|\langle F, G \rangle_{\widetilde L}|\leq 
\|F\|_{\widetilde W^{a, p}_q} \|G\|_{\widetilde W^{-a, p'}_q}
\end{equation}
see again \cite{A}. Note that since
\begin{equation*}
 \widetilde C^{\infty, c}_q(\mathbb R^N) \ \cap\ 
 C^{\infty, c}_q(\mathbb R^N) = \{0\}
\end{equation*}
the spaces $\widetilde L^r_q(\mathbb R^N)$  and
$L^r_q(\mathbb R^N)$ are transversal;
  the same is true for $\widetilde W^{a, p}_q(\mathbb R^N)$ and 
  $ W^{a,p}_q(\mathbb R^N)$
 and for the respective dual spaces.
\medskip

\subsection{Operators and their adjoints}\label{Ss:operators}
For $\aleph$ as in \eqref{E:correspond} and for any admissible degree increment $\ell$, 
we define 
   a $k$th-order differential condition via the action
   \begin{equation}\label{E:gen-LR-oper}
 \sum\limits_{I\in\mathcal I(N, q)}\!\!\!\!F_I\,dz^I\quad
  \to \ \ 
 \sum\limits_{L\in\mathcal I(N, \,q+\ell)}\!\!\bigg(
  \sum\limits_{\atopp{I\in\mathcal I(N, q)}{\alpha\in\mathcal S(n, k)}}
  \!\!\!
  \epsilon_L^{\aleph(i\alpha) I}\
  \frac{\dee^kF_I}{\dee z^{i\alpha}}
   \bigg)dz^L
   \end{equation}
  where $N$ is as in \eqref{E:embedding-dim} and 
 \eqref{E:constraint-dim} and $ q\in\{0, 1, \ldots, N\}$. Here 
 $i\alpha$
 is as in \eqref{E:multi-index-embed-2}
 and
 $\aleph$ is the correspondence \eqref{E:correspond}.
   This action produces a differential operator $T_{\ell, \aleph}$ that maps
      \begin{equation}\label{E:T-q-map}
 T_{\ell,\, \aleph}:  C^{\infty, c}_q(\mathbb R^N)\to C^{\infty, c}_{q+\ell}(\mathbb R^N), \quad 0\leq q\leq N.
   \end{equation}
It follows from \eqref{E:sobolev-deriv-3} that the action \eqref{E:gen-LR-oper} also
determines an operator
    \begin{equation}\label{E:tilde-T-q-map}
 \widetilde T_{\ell,\, \aleph}:  \widetilde C^{\infty, c}_q(\mathbb R^N)\to
  \widetilde C^{\infty, c}_{q+\ell}(\mathbb R^N), \quad 0\leq q\leq N.
   \end{equation}
 Now observe that \eqref{E:pi-i} grants   that the pullback by $\pi$ maps
   $$\pi^*: C^{\infty, c}_q(\mathbb R^n)\to
   \widetilde C^{\infty, c}_q(\mathbb R^N)$$
 see \eqref{E:canonical-proj}. On the other other hand, it is immediate to check that
  $$
  i^*: \widetilde C^{\infty, c}_q(\mathbb R^N)\to C^{\infty, c}_q(\mathbb R^n).
  $$
On account of these observations we see that the action 
\eqref{E:gen-LR-oper} produces a third operator 
$\Top$ that maps 
 \begin{equation}\label{E:Top-map}
\Top:  C^{\infty, c}_q(\mathbb R^n)\to
  C^{\infty, c}_{q+\ell}(\mathbb R^n), \quad 0\leq q\leq n
   \end{equation}
  and is defined as follows
\begin{equation}\label{E:Top-def}
\Top := i^*\,\widetilde T_{\ell,\aleph}\, \pi^*\, .
\end{equation}
    Note that $\Top$ acts non-trivially only for
 \begin{equation}\label{E:constr-n-ell}
 n\geq \ell.
 \end{equation}
Condition \eqref {E:constr-n-ell} may be viewed in two different ways: 
as a constraint on the size of the degree increment $\ell$ relative to 
the pair $(n, k)$
(however note that \eqref{E:constr-n-ell} 
is satisfied by $\ell =1$ for {\it any} pair $(n, k)$)
or as a constraint on the size of $n$ relative to 
$k$
 (and in this case, imposing the constraint
$n\geq k$ ensures that \eqref{E:constr-n-ell} holds for {\em all} admissible degree increments).

\smallskip

\noindent In the following, $\langle \cdot, \cdot\rangle$ denotes the duality
in $W^{p, 2}_q(\mathbb R^N)$ (resp. $W^{p, 2}_q(\mathbb R^n))$.
  \begin{prop}\label{P:hilbert-T-adj} Let $\ell$ be an admissible degree increment for $(n, k)$.
  
The formal adjoint of $T_{\ell, \aleph}$ on $W^{a,2}_q(\mathbb R^N)$ is
\begin{equation}\label{E:def-hilbert-T-adj}
T^*_{\ell,\aleph} := (-1)^{k+q(N-\ell-q)}*_NT_{\ell, \aleph}\,*_N,\quad 
0\leq q\leq N
\end{equation}
 That is, for any $F\in C^{\infty, c}_q(\mathbb R^N)$ and for any
 $G\in C^{\infty, c}_{q+\ell}(\mathbb R^N)$ we have
\begin{equation}\label{E:T-hilbert-adj}
\langle\, T_{\ell, \aleph}F,\, G\,\rangle \ =\ 
\langle\, F,\, T^*_{\ell,\aleph} G\,\rangle\, 
\end{equation}
The formal adjoint of $\widetilde T_{\ell, \aleph}$ on $\widetilde W^{a, 2}_q(\mathbb R^N)$ is 
 \begin{equation}\label{E:def-hilbert-tilde-T-adj}
\widetilde T^*_{\ell,\aleph} :=
(-1)^{k+q(N-\ell-q)}*_N\widetilde T_{\ell, \aleph}\,*_N,
\quad 
0\leq q\leq N
\end{equation}
 That is,  for any $F\in \widetilde C^{\infty, c}_q(\mathbb R^N)$ and for any
 $G\in \widetilde C^{\infty, c}_{q+\ell}(\mathbb R^N)$ we have
\begin{equation}\label{E:hilbert-tilde-T-adj}
\langle\, \widetilde T_{\ell, \aleph}F,\, G\,\rangle_{\widetilde L} \ =\ 
\langle\, F,\, \widetilde T^*_{\ell,\aleph} G\,\rangle_{\widetilde L}\, 
\end{equation}
Suppose further that $n\geq \ell$. Then, the formal adjoint of $\Top$ on 
$W^{a, 2}_q(\mathbb R^n)$ is 
 \begin{equation}\label{E:def-Top-adj}
\AdjTop := (-1)^{k+q(n-\ell -q)}*_n\Top*_n
\quad 
0\leq q\leq n
\end{equation}
that is,  for any $f\in C^{\infty, c}_q(\mathbb R^n)$ and for any
 $g\in  C^{\infty, c}_{q+\ell}(\mathbb R^n)$ we have
\begin{equation}\label{E:hilbert-Top-adj}
\langle\,\Top f,\, g\,\rangle \ =\ 
\langle\, f,\, \AdjTop g\,\rangle\, .
\end{equation}
 \end{prop}
 \bigskip

\subsection{Estimates}\label{Ss:estimates} 
\medskip

\begin{thm}\label{P:LR-inflated-dim}
   Let $n\ge 2$ and $k\ge 1$ be given. Let $\ell\in\{1,\ldots, k\}$ be any admissible degree increment for the pair $(n, k)$, and 
   let $N$ be a solution of 
   \eqref{E:embedding-dim} that satisfies \eqref{E:constraint-dim}.  
     \begin{itemize}
  \item[] For any $0\leq q\leq N-\ell$ and for any 
  $F\in\widetilde L^1_q(\mathbb R^N)$, 
   if
   \begin{equation}\label{E:LR-div-invar-k-inflated}
   \widetilde T_{\ell,\, \aleph}\,F = 0
   \end{equation}
   in the sense of distributions, then 
   \begin{equation}\label{E:LR-ineq-inflated}
   \left|\langle F, H\rangle_{\widetilde L}\right|
  \leq 
  \widetilde C\| F\|_{\widetilde L^1_q}\,\|\nabla H\|_{\widetilde L^n_q}
    \end{equation}
   for any 
   $H\in (\widetilde L^\infty_q \cap \widetilde W^{1, n}_q)(\mathbb R^N)$.
   \bigskip
   
     \item[] For any $\ell\leq p\leq N$ and for any
     $G\in\widetilde L^1_p(\mathbb R^N)$, 
 if
   \begin{equation}\label{E:tilde-T-star-LR}
  \widetilde T_{\ell,\, \aleph}^* G = 0
   \end{equation}
   in the sense of distributions, then 
   \begin{equation}\label{E:LR-dual-ineq-inflated}
     \left|\langle G, K\rangle_{\widetilde L}\right|
   \leq 
  \widetilde C\| G\|_{\widetilde L^1_p}\,\|\nabla K\|_{\widetilde L^n_p}
   \end{equation}
    for any 
   $K\in (\widetilde L^\infty_p \cap \widetilde W^{1, n}_p)(\mathbb R^N)$.
   \end{itemize}
   \medskip
   
      The constant $\widetilde C$ depends only on $n$ and $k$.
   \end{thm}
   \medskip

\begin{thm}\label{P:LR-source-smallest}
 Let $n\ge 2$ and $k\ge 1$ be given. Let $\ell\in\{1,\ldots, k\}$ be an admissible degree increment for the pair $(n, k)$ such that
 $$n\geq \ell\, .$$ 
  \begin{itemize}
  \item[] 
  For any $0\leq q\leq n-\ell$ and for any 
  $ f\in  L^1_q(\mathbb R^n)$, if
 \begin{equation}\label{E:LR-div-invar-k-source-smallest}
   \Top f = 0
   \end{equation}
   in the sense of distributions, then 
   \begin{equation}\label{E:LR-ineq-source-smallest}
   \left|\langle f, h\rangle\right|\leq 
   \mathcal C 
      \|f\|_{L^1_q}\,\|\nabla h\|_{L^n_q}
   \end{equation}
   for any $h\in (L^\infty_q \cap W^{1, n}_q)(\mathbb R^n)$.
   \bigskip
   
   \item[]
    For any $\ell\leq p\leq n$ and for any 
  $ g\in  L^1_p(\mathbb R^n)$, if
 \begin{equation}\label{E:LR-div-invar-k-source-smallest-adj}
   \AdjTop g = 0
   \end{equation}
   in the sense of distributions, then 
   \begin{equation}\label{E:LR-ineq-source-smallest-adj}
   \left|\langle g, h\rangle\right|\leq 
   \mathcal C
    \|g\|_{L^1_p}\,\|\nabla h\|_{L^n_p}
   \end{equation}
   for any $h\in (L^\infty_p \cap W^{1, n}_p)(\mathbb R^n)$.
   \end{itemize}
   \medskip
   
   The constant $\mathcal C$ depends only on $n$ and $k$.
\end{thm}
\medskip

We have
\bigskip

\centerline{Theorem \ref{T:VS-thm-4}  
 $\iff$
Theorem \ref{P:LR-inflated-dim}
$\Rightarrow$
Theorem \ref{P:LR-source-smallest}
$\Rightarrow$
Theorem \ref{T:BB-thm-3}.}
\bigskip

\subsection{Hodge systems}\label{Ss:hodge}
Concerning the compatibility conditions for the Hodge system
for each of $\widetilde T_{\ell, \aleph}$ and $\Top$, we have
\begin{equation*}
(\widetilde T_{\ell, \aleph}\circ\widetilde T_{\ell, \aleph})\, F
= (1+ (-1)^{\ell^2})
  \!\!\!  \!\!\!\!
\sum\limits_{M\in\mathcal I(N, q+2\ell)}\!\!\!
\left(
\sum\limits_{\atopp{I\in\mathcal I(N, q)}{\alpha, \beta \in\mathcal S(n, k)}}
  \!\!\!
  \epsilon^{\aleph(i\beta)\aleph(i\alpha)I}_M
  \frac{\dee^{2k}F_I}{\dee z^{i\alpha}\dee z^{i\beta}}
  \right)\!dz^M
\end{equation*}
so in particular
\begin{equation}\label{E:compatibility-big-T}
\widetilde T_{\ell, \aleph}\circ\widetilde T_{\ell, \aleph}\,
=\,0 \iff \ell\ \ \mbox{is odd}.
\end{equation}
A similar computation shows that the same
 is true for $\Top$, so in the sequel we will often pay special attention to the admissible degree increment
$\ell = 1$.
\begin{lem}\label{L:box-reps}
Let $\widetilde T^*_{\ell,\aleph}$ be given by \eqref{E:def-hilbert-tilde-T-adj}
and set
\begin{equation}\label{E:big-box-def}
\widetilde\square_{\ell, \aleph} := \widetilde T_{\ell, \aleph}\,
\widetilde T^*_{\ell, \aleph} + 
\widetilde T^*_{\ell, \aleph}\,\widetilde T_{\ell, \aleph}\ .
\end{equation}
If 
$$
H=\sum\limits_{I\in\mathcal I(N, q)}\!\!\!\! H_I\,dz^I\in \widetilde C^{\infty, c}_q(\mathbb R^N)
$$
then
\begin{equation}\label{E:big-box-ell-repr}
\widetilde\square_{\ell, \aleph} H =\
(-1)^{k  + \ell N}\!\!\!\!\!\!
\sum\limits_{\atopp{M, I\in\mathcal I(N, q)}{\alpha, \beta \in\mathcal S(n, k)}}
\!\!\!\!\!\widetilde C^{MI}_{\aleph(i\alpha)\,\aleph(i\beta)}\,
\frac{\dee^{2k}H_I}{\dee z^{i\alpha}\dee z^{i\beta}}\, dz^M\, 
\end{equation}
where
\begin{equation}\label{E:big-coeffs}
\widetilde C^{MI}_{\aleph(i\alpha)\,\aleph(i\beta)} =
\sum\limits_{L\in\mathcal I(N, q+\ell)}
\!\!\!\!\!\! \epsilon^L_{\aleph(i\alpha)I}\!\cdot\epsilon^L_{\aleph(i\beta) M}
\ +
\!\!\!\sum\limits_{K\in\mathcal I(N, q-\ell)}
\!\!\!\!\!\! \epsilon^M_{\aleph(i\alpha) K}\!\cdot\epsilon^I_{\aleph(i\beta) K}.
\end{equation}
In particular, for $\ell =1$ we have
\begin{equation}\label{E:big-box-1-repr}
\widetilde\square_{1, \aleph} H =
\sum\limits_{I\in\mathcal I(N, q)}\!\!\!\!
(\widetilde\square_{1, \aleph}H_I)\,dz^I = (-1)^{k  + N}\!\!\!
\sum\limits_{\atopp{I\in\mathcal I(N, q)}{\alpha\in\mathcal S(n, k)}}
  \!\! \frac{\dee^{2k}H_I}{\dee z^{i\alpha}\dee z^{i\alpha}}\,dz^I\, .
\end{equation}
Let $\Top$ be given by \eqref{E:def-Top-adj} (assume $n\geq \ell$) and set
\begin{equation}\label{E:box-def}
\square_{\ell,\aleph} :=
\Top\Top^* + \Top^*\Top.
\end{equation}
If 
$$
h=\sum\limits_{I\in\mathcal I(n, q)}\!\!\!\! h_I\,dz^I\in  C^{\infty, c}_q(\mathbb R^n)
$$
then
\begin{equation}\label{E:box-ell-repr}
\square_{\ell, \aleph} h =\
(-1)^{k  + \ell n}\!\!\!\!\!\!\!\!\!\!\!\!
\sum\limits_{\atopp{M, I\in\mathcal I(n, q)}
{\alpha, \beta\, \in\, (\pi\circ\aleph^{-1})(\mathcal I(n, \ell))}}
\!\!\!\!\!\!\!\!\!\!\!\!\!C^{MI}_{\aleph(i\alpha)\,\aleph(i\beta)}\,
\frac{\dee^{2k}h_I}{\dee x^\alpha\dee x^\beta}\, dx^M\, 
\end{equation}
where
\begin{equation}\label{E:coeffs}
C^{MI}_{\aleph(i\alpha)\,\aleph(i\beta)} =
\sum\limits_{L\in\mathcal I(n, q+\ell)}
\!\!\!\!\!\! \epsilon^L_{\aleph(i\alpha)I}\!\cdot\epsilon^L_{\aleph(i\beta) M}
\ +
\!\!\!\sum\limits_{K\in\mathcal I(n, q-\ell)}
\!\!\!\!\!\! \epsilon^M_{\aleph(i\alpha) K}\!\cdot\epsilon^I_{\aleph(i\beta) K}.
\end{equation}
In particular, for $\ell =1$ we have
\begin{equation}\label{E:box-1-repr}
\square_{1, \aleph} h\ =
\sum\limits_{I\in\mathcal I(n, q)}\!\!\!\!
(\square_{1, \aleph}h_I)\,dx^I =\ (-1)^{k  + n}\!\!\!\!\!\!\!\!\!\!\!\! 
\sum\limits_{\atopp{I\in\mathcal I(n, q)}
{\alpha \in\, (\pi\circ\aleph^{-1})(\{1,\ldots, n\})}}\!\!\!\!\!\!\!\!\! 
  \!\! \frac{\dee^{2k}h_I}{\dee x^\alpha \dee x^\alpha}\,dx^I\, .
\end{equation}
\end{lem}

\begin{cor}\label{C:Legendre-Hadamard} For any $0\leq q\leq N$ and for any choice of the correpondence $\aleph$, the operator 
$\widetilde\square_{1,\aleph}$ satisfies the Legendre-Hadamard condition  in the following sense:
\begin{equation}\label{E:big-box-elliptic}
\mathrm{Re}\left(\sum\limits_{\atopp{I, M\in\mathcal I(N, q)}{\alpha, \beta\in\mathcal S(n, k)}}
  \!\!\!\!\!\! \widetilde C^{MI}_{\aleph(i\alpha)\aleph(i\beta)}\,\xi^{i\alpha}\xi^{i\beta}
  \zeta_I\overline{\zeta}_M\right)\geq \ C\,|\xi|^{2k}|\zeta|^2
\end{equation}
for any $\xi\in i(\mathbb R^n)$ and for any 
$\displaystyle{\zeta\in\mathbb C^{\binom{N}{q}}}$.
\end{cor}
See \cite{KPV}. Indeed, by \eqref{E:big-box-1-repr} we have that the coordinate-based 
representation of $\widetilde\square_{1,\aleph}$ is independent of the choice of $\aleph$ and furthermore
$$
\mbox{Re}\left(\sum\limits_{\atopp{I, M\in\mathcal I(N, q)}{\alpha, \beta\in\mathcal S(n, k)}}
  \!\!\!\!\!\! \widetilde C^{MI}_{\aleph(i\alpha)\aleph(i\beta)}\,\xi^{i\alpha}\xi^{i\beta}
  \zeta_I\overline{\zeta}_M\right) = 
  \left(\sum\limits_{\alpha\in\mathcal S(n, k)}
   \!\!\!\! \xi_1^{2\alpha_1}\cdots\xi_n^{2\alpha_n}
    \right)|\zeta|^2\
$$
and if $\xi\in i(\mathbb R^n)$ then 
$$
\sum\limits_{\alpha\in\mathcal S(n, k)}
   \!\!\!\! \xi_1^{2\alpha_1}\cdots\xi_n^{2\alpha_n}
\geq C|\xi|^{2k}
$$
where $C=C(n, k)$.
On the other hand, the coordinate-based representation of
$\square_{1,\aleph}$ does depend on the choice of $\aleph$, see 
\eqref{E:big-box-1-repr}, and so does the
uniform ellipticity of $\square_{1,\aleph}$; for instance, if $\aleph$ is chosen 
so that 
$$
(\pi\circ\aleph^{-1})(\{1, 2,\ldots, n\}) =$$
$$
 =\{(k,0,\ldots,0),(1,k-1,0,\ldots, 0),\ldots, (1,0,\ldots,0, k-1)\}
\subset\mathcal S(n, k)
$$
then $\square_{1,\aleph}$ has symbol
$$
\left(\sum\limits_{j=1}^n\xi_1^2\,\xi_j^{2(k-1)}\right)|\zeta|^2,\quad 
\xi\in\mathbb R^n,\ \eta\in\mathbb C^{\binom{n}{q}}
$$
which fails to be uniformly elliptic (take e.g., $\xi :=(0, 1,\ldots, 1)$). Choosing
instead $$
(\pi\circ\aleph^{-1})(\{1, 2,\ldots, n\}) =$$
$$
\{(k, 0,\ldots,0), (0, k,0,\ldots, 0),\ldots, (0,\ldots, 0, k)\}
$$
(corresponding to the example $\mathcal T$ discussed in the Introduction)
leads to an operator $\square_{1,\aleph}$ which is easily verified to be uniformly elliptic, as we have
$$
\sum\limits_{j=1}^n\xi_j^{2k}\ge n^{1-k}|\xi|^{2k}.
$$
\begin{lem}\label{L:big-box-elliptic-estimates}
We have that
\begin{equation*}
\widetilde\square_{1,\aleph}\!\!: \widetilde C^{\infty, c}_q(\mathbb R^N)\to  \widetilde C^{\infty, c}_q(\mathbb R^N)
\end{equation*}
is invertible
for any 
$1<p<\infty$. 
For any $\varphi\in
\widetilde L^p_q (\mathbb R^N)$  we have
\begin{equation}\label{E:big-box-elliptic estimates}
\|\Phi\|_{\widetilde W_q^{2k, p}}\lesssim
\|\varphi\|_{\widetilde L_q^{ p}}
\end{equation}
where $\Phi := \widetilde\square^{-1}_{1, \aleph}\,\varphi$.
\end{lem}
\medskip

\begin{thm}\label{T:LR-thm-Z} 
 Suppose that $F\in \widetilde L^1_{q+1}(\mathbb R^N)$ and
   $G \in \widetilde L^1_{q-1}(\mathbb R^N)$ 
   satisfy the hypotheses of Theorem \ref{P:LR-inflated-dim}. Let
 \begin{equation}
 Z = \widetilde \square_{1,\aleph}^{-1}(\widetilde T_{1,\aleph}^*F + \widetilde T_{1,\aleph}G)\in\widetilde\Lambda_q(\mathbb R^N),\ \ 0\leq q\leq N
 \end{equation}
 be the solution of the Hodge system
  for $\widetilde T_{1,\aleph}$ with data $(F, G)$, that is:
 \begin{equation}
 \label{E:F-q}
 \bigg\{
 \begin{array} {ll}
 \widetilde T_{1,\aleph} Z \ =& F\\
 \widetilde T_{1,\aleph}^*Z\ =& G
 \end{array}
 \end{equation}
 Then 
   \begin{equation}
 \label{E:ineq-F-q}
 \|Z\|_{\widetilde W^{k-1, r}_q}\leq C\big(\|F\|_{\widetilde L^1_{q+1}} + 
 \|G\|_{\widetilde L^1_{q-1}}\big),\ \ \mbox{for} \ \  r= n/(n-1)
  \end{equation}
 whenever $q$ is neither 1 (unless $G=0$) nor $N-1$ (unless $F=0$).
 \end{thm}
\bigskip

We have:
\bigskip

\centerline{Theorem \ref{P:LR-inflated-dim} ($\ell =1$;\quad  $1\leq p,\, q\leq N-1$)\ $\iff$\ 
Theorem \ref{T:LR-thm-Z}}
\bigskip

For those choices of $\aleph$ that give rise to a uniformly elliptic
 $\square_{1,\aleph}$, an analogous result holds for
 \begin{equation}\label{E:f-q}
\left\{ \begin{array}{rcl}
\mathcal T_{1,\aleph} h =&f,  &\mathcal T_{1,\aleph}f =0\\
 \mathcal T_{1,\aleph}^*h =&g, &\mathcal T_{1,\aleph}^*g=0
\end{array}\right.
\end{equation}
which turns out to be equivalent to Theorem \ref{P:LR-source-smallest} ($\ell =1$). We omit the statement.

\bigskip

We remark in closing that for $\ell\geq 2$ there is no analog of \eqref{E:big-box-1-repr}. Indeed, setting
$$
\{\lambda_0\} := \{\aleph(i\alpha)\}\cap \{\aleph(i\beta)\}
$$
and
$$
\{\widehat{\aleph(i\alpha)}\} := 
\{\aleph(i\alpha)\} \setminus \{\lambda_0\}\,;\quad 
\{\widehat{\aleph(i\beta)}\}:=  \{\aleph(i\beta)\} \setminus \{\lambda_0\}
$$
(where the brackets $\{\ \}$ indicate that the (ordered) label $J$ is being regarded as an (unordered) set $\{J\}$), it can be proved that
\begin{equation}
\widetilde C^{MI}_{\aleph(i\alpha) \aleph(i\beta)} =
(1+ (-1)^{(\ell - |\{\lambda_0\}|)^2})\,\epsilon^
{\aleph(i\alpha)}_{\lambda_0\widehat{\aleph(i\alpha)}}\,
\cdot\epsilon^{\aleph(i\beta)}_{\lambda_0\widehat{\aleph(i\beta)}}\,\cdot\epsilon^{\widehat{\aleph(i\alpha)} I}_{\widehat{\aleph(i\beta)} M}.
\end{equation} 
In particular, the coordinate-based representation of $\widetilde\square_{\ell, \aleph}$ does depend on the choice of the representation $\aleph$, and it is no longer true that  $\widetilde C^{MI}_{\aleph(i\alpha) \aleph(i\beta)}= 0$ whenever $\alpha\neq \beta$, even for odd $\ell$. 
 \bigskip

\section{Proofs}\label{S:Proofs}
\medskip

{\it Proof of Lemma \ref{L:hilbert}}.\quad 
Conclusions {\em i.} and {\em ii.} are an immediate consequence of the (classical) theory for $\mathbb R^n$ combined with the
 readily verified identities:
\begin{equation}\label{E:r-big-norm-1}
\|F\|_{\widetilde L^r_q}^{^r} =
\sum\limits_{I\in\mathcal I(N, q)}\
\int\limits_{\mathbb R^n}\!\!
\bigg(\sum\limits_{I\in\mathcal I(N, q)}|F_I\circ\,i|^2(x)\bigg)^{\! r/2}\!dV(x)
\end{equation}
and
\begin{equation}\label{E:pairing-big-1}
\langle F, G\rangle_{\widetilde L} =
\sum\limits_{I\in\mathcal I(N, q)}\
\int\limits_{\mathbb R^n}
(F_I\circ i)(x)\cdot (G_I\circ i)(x)\,dV(x)\,.
\end{equation}
To prove the density of $\widetilde C^{\infty, c}_q(\mathbb R^N)$ in
 $\widetilde L^r_q(\mathbb R^N)$, let 
 $$F= \sum\limits_{I\in\mathcal I(N, q)} F_I\, dz^I\in \widetilde L^r_q(\mathbb R^N)$$
be given. By the definition of $\widetilde L^r_q(\mathbb R^N)$, for any $I\in\mathcal I(N, q)$ we have that $F_I\circ i\in L^r(\mathbb R^n)$, and so  there is 
$\{f_{j, I}\}_{j\in\mathbb N}\subset C^\infty_c(\mathbb R^n)$
such that
\begin{equation}\label{E:density-1}
\|\, f_{j, I}- F_I\circ i\,\|_{L^r(\mathbb R^n)}\to 0\quad \mbox{as}\ j\to\infty.
\end{equation}
Define
$$
F_j=
\sum\limits_{I\in\mathcal I(N, q)}
\!\!\!
F_{j, I}\, dz^I,\quad F_{j, I}:= f_{j, I}\circ\pi
$$
Then, using \eqref{E:pi-i}, we see that
$$
F_{j, I}\circ\, i\circ\pi = F_{j, I}\quad \mbox{and}\quad 
F_{j, I}\circ\, i = f_{j, I}\in C^\infty_c(\mathbb R^n)
$$
hold for any $I\in\mathcal I(N, q)$,
and from these it follows that
$$
\{F_j\}_{j\in\mathbb N}\subset\widetilde C^{\infty, c}_q(\mathbb R^N).
$$
Moreover, on account of  \eqref{E:r-big-norm-1} and 
\eqref{E:density-1},
there is $C=C(r, N)$ such that
$$
\| F_j - F\|^r_{\widetilde L^r_q}\ \leq \ C \!\!\!\!
\sum\limits_{I\in\mathcal I(N, q)}
\!\!\!
\| f_{j, I} - F_I\circ\, i\ \|^r_{L^r(\mathbb R^n)}\to 0\quad \mbox{as}\ \ j\to\infty,
$$
 as desired. 
The conclusions concerning the Sobolev spaces
follow from the theory for
 $W^{a, p}_q(\mathbb R^n)$ 
  via 
\eqref{E:sobolev-deriv-2}.
\qed
 \bigskip

{\it Proof of  Proposition \ref{P:hilbert-T-adj}}.\quad Let $F\in C^{\infty, c}_q(\mathbb R^N)$ and $G\in C^{\infty, c}_{q+\ell}$ be given.
One has
$$
T_{\ell, \aleph} F\wedge *_N G\ \,=
\sum\limits_{L\in \mathcal I(N, q+\ell)}
\bigg(\sum\limits_
{\atopp{I\in\mathcal I(N, q)}{\alpha\in\mathcal S(n, k)}}\!\!
\epsilon^{\aleph(i\alpha)I}_L
\frac{\dee^k F_I}{\dee z^{i\alpha}}\,G_L\bigg)dV
$$
Integrating
both sides of this identity and then further integrating the right-hand side by parts $k$-many times we find
\begin{equation}\label{E:auxx}
\langle\, T_{\ell, \aleph}F, G\,\rangle = (-1)^k\int\limits_{\mathbb R^N}
\sum\limits_{L\in \mathcal I(N, q+\ell)}
\bigg(\sum\limits_
{\atopp{I\in\mathcal I(N, q)}{\alpha\in\mathcal S(n, k)}}\!\!
\epsilon^{\aleph(i\alpha)I}_L F_I\,
\frac{\dee^k G_L}{\dee z^{i\alpha}}\bigg)dV
\end{equation}
On the other hand, a computation that requires manipulating the 
coefficients $\epsilon^{\aleph(i\alpha)K}_J$ gives
\begin{equation*}
F\wedge *_N\big(*_N T_{\ell, \aleph} *_N\ G\big) =
(-1)^{q(N-q) + q\ell}\!\!\!\!\!\!\!\!
\sum\limits_{L\in \mathcal I(N, q+\ell)}\!\!\!
\bigg(\sum\limits_
{\atopp{I\in\mathcal I(N, q)}{\alpha\in\mathcal S(n, k)}}\!\!
\epsilon^{\aleph(i\alpha)I}_L F_I\,
\frac{\dee^k G_L}{\dee z^{i\alpha}}\bigg)dV
\end{equation*}
Identity \eqref{E:T-hilbert-adj} is now obtained by integrating the two sides of the identity 
above and comparing with
\eqref{E:auxx} after having adjusted the multiplicative constants as in \eqref{E:def-hilbert-T-adj}. 
Note that since 
$$
D^\lambda\, T_{\ell, \aleph} F = T_{\ell, \aleph}\, D^\lambda F \quad \mbox{for any multi-index}\ \lambda
$$
where $D^\lambda F\in C^{\infty, c}_q(\mathbb R^N)$ is defined as in 
\eqref{E:def-deriv-form}, the same argument also shows that
$$
\langle D^\lambda\, T_{\ell, \aleph} F,\, D^\lambda G\rangle =
\langle D^\lambda F,\, D^\lambda\, T^*_{\ell, \aleph} G\rangle
$$
The proofs of  \eqref{E:hilbert-tilde-T-adj} and of \eqref{E:hilbert-Top-adj} follow in a similar fashion.
\qed

\bigskip

   {\it {Theorem \ref{T:VS-thm-4}}} $\Rightarrow$
           {\it {Theorem \ref{P:LR-inflated-dim}}}.\quad 
          Let $\ell$ be an admissible degree increment and let 
          $0 \leq q\leq N-\ell$. Suppose that 
    $$F=\sum\limits_{I\in\mathcal I(N, q)}\!\!\! F_I\,dz^I\in \widetilde\Lambda_q(\mathbb R^N)$$ 
    and
        $$H=\sum\limits_{I\in\mathcal I(N, q)}\!\!\!H_I\,dz^I\in\widetilde\Lambda_q(\mathbb R^N)$$ 
     satisfy the hypotheses of Theorem \ref{P:LR-inflated-dim}.
        Fix an arbitrary $I_0\in\mathcal I(N, q)$, and 
       choose (any) $L_0\in\mathcal I(N, q+\ell)$
       so that
    \begin{equation*}
    I_0\subseteq L_0.
    \end{equation*}
    (The hypothesis: $q\leq N-\ell$ 
        grants $q+\ell \leq N$ and so at least one such $L_0$ must exist.)
    With $I_0$ and $L_0$ fixed as above, define
    $\displaystyle{h^{L_0}= (h^{L_0}_\alpha)_{\alpha\in\mathcal S(n, k)}}$ 
    and $\displaystyle{g^{L_0}= (g^{L_0}_\alpha)_{\alpha\in\mathcal S(n, k)}}$ via
     \begin{equation*}
    h_\alpha^{L_0} =\sum\limits_{I\in\mathcal I(N, q)}\!
    \epsilon^{\aleph(i\alpha) I}_{L_0}\,H_I\circ i,\qquad \alpha \in \mathcal S(n, k);
    \end{equation*}
    \begin{equation*}
    g_\alpha^{L_0} =\sum\limits_{I\in\mathcal I(N, q)}\!
    \epsilon^{\aleph(i\alpha)I}_{L_0}\, F_I\circ i,\qquad \alpha \in \mathcal S(n, k).
    \end{equation*}
     
    We claim that $g^{L_0}$ satisfies condition \eqref{E:VS-div-k}
     in Theorem \ref{T:VS-thm-4}: to this end, note that by 
     \eqref{E:sobolev-deriv-3} we have
          $$
           \sum\limits_{\alpha\in \mathcal S(n, k)}\frac{\dee^k g^{L_0}_\alpha}{\dee x^\alpha}=
  \left(\sum\limits_{\atopp{I\in\mathcal I(N, q)}{\alpha\in\mathcal S(n, k)}}
  \!\!\!\!\!
  \epsilon_{L_0}^{\aleph(i\alpha) I}\,
  \frac{\dee^kF_I}{\dee z^{i\alpha}}\right)\circ i = [T_{\ell,\,\aleph} F]_{L_0}\circ i = 0
      $$
where the last identity is due to the
     hypothesis \eqref{E:LR-div-invar-k-inflated}.
    It thus follow from Theorem \ref{T:VS-thm-4} that
    \begin{equation*}
    \bigg|\,\int\limits_{\mathbb R^n}g_{\alpha_0}^{L_0}\cdot h_{\alpha_0}^{L_0}
    \,\bigg|\leq C\|g^{L_0}\|_{L^1(\mathbb R^n)}\|\nabla h^{L_0}\|_{L^n(\mathbb R^n)}
    \end{equation*}
where $\alpha_0\in\mathcal S (n, k)$ is uniquely determined by $I_0$ and $L_0$ via
$$
i\alpha_0 := \aleph^{-1}(L_0\setminus I_0),
$$
(note that 
  $
 L_0\setminus I_0\in\mathcal I(N, \ell)$.)
But for $\alpha_0$ as above we have
$$
\epsilon^{\aleph (i\alpha_0)I}_{L_0}\neq 0 \iff I =I_0
$$
and so 
\begin{equation*}
g_{\alpha_0}^{L_0} = \epsilon^{\aleph(i\alpha_0)I_0}_{L_0}\,F_{I_0}\circ i,
\quad \text{and}\quad
h_{\alpha_0}^{L_0} = \epsilon^{\aleph(i\alpha_0)I_0}_{L_0}\,H_{I_0}\circ i.
\end{equation*}
On the other hand, it is immediate to verify that
\begin{equation*}
\|g^{L_0}\|_{L^1(\mathbb R^n)}\lesssim \|F\|_{\widetilde L^1_q(\mathbb R^N)},\quad
\|\nabla h^{L_0}\|_{L^n(\mathbb R^n)}\lesssim 
\|\nabla H\|_{\widetilde L^n_q(\mathbb R^N)},
\end{equation*}
and
$$
\int\limits_{\mathbb R^n}g_{\alpha_0}^{L_0}\cdot h_{\alpha_0}^{L_0} = 
\int\limits_{\mathbb R^n} (F_{I_0}\circ i)\cdot (H_{I_0}\circ i)(x) dV(x).
$$
Since $I_0\in \mathcal I(N, q)$ had been fixed arbitrarily, we have proved that
\begin{equation}\label{E:aux1}
\left|\, \int\limits_{\mathbb R^n} \Big((F_{I}\circ i)\cdot (H_{I}\circ i)\Big)(x) dV(x)\right|
\leq C \|F\|_{\widetilde L^1_q(\mathbb R^N)}\, \|\nabla H\|_{\widetilde L^n_q(\mathbb R^N)}
\end{equation}
holds for any $I\in\mathcal I(N, q)$, for any $0\leq q\leq N-\ell$ and for any admissible degree increment $\ell$. Inequality \eqref{E:LR-ineq-inflated} follows from  \eqref{E:aux1}
and the coordinate-based representation for $\langle\cdot, \cdot\rangle_{\widetilde L}$, see \eqref{E:pairing-big-1}.
(We remark that in the special case $q=0$, the proof follows along these very same lines by defining $g_\alpha := F\circ i$ for each $\alpha \in \mathcal S(n, k)$.)
\smallskip

In order to prove \eqref{E:LR-dual-ineq-inflated}, it suffices to apply
\eqref{E:LR-ineq-inflated} to: $F:= *_NG$ and $H:=*_N K$ (with $q:=N-p$).
\qed

\bigskip

{\it Theorem \ref{P:LR-inflated-dim} $\Rightarrow$ Theorem \ref{T:VS-thm-4}.}\quad 
 Let $\ell$ be any admissible
degree increment for $(n, k)$ and let $\aleph$ be any one-to-one correspondence: $i\mathcal S (n, k)\to \mathcal I(N,\ell)$. Suppose that $g$ and $h$ satisfy the hypotheses of 
 Theorem \ref{T:VS-thm-4}; without loss of generality we may assume that
 $g_\alpha, \ h_\alpha \in C^{\infty}_0(\mathbb R^n)$, $\alpha \in\mathcal S(n, k)$. 
Choose $q:= N-\ell$ and define $F$ and $H$ in $ \Lambda_{N-\ell}(\mathbb R^N)$ via
\begin{equation}\label{E:def-F_I}
F_I:= \epsilon_{(1,\ldots, N)}^{I'I}\,g_\alpha\,\circ\, \pi,
\quad I\in \mathcal I(N, N-\ell)
\end{equation}
\begin{equation}\label{E:def-H_I}
H_I:= \epsilon_{(1,\ldots, N)}^{I'I}\,h_\alpha\,\circ\, \pi,
\quad I\in \mathcal I(N, N-\ell)
\end{equation}
where
$I':= \{1,\ldots, N\}\setminus I\in\mathcal I(N,\ell)$, and $\alpha\in\mathcal S(n, k)$ is uniquely determined
 by $I$ and by $\aleph$ via
$$
i\alpha =\aleph^{-1}(I').
$$
Since $\pi\,\circ\, i$ is the identity on $\mathbb R^n$, we have
$$
F_I\,\circ\, i\, \circ\, \pi = F_I,\quad F_I\circ i = \epsilon^{I' I}_{(1,\ldots, N)}\,g_\alpha\in C^{\infty}_0(\mathbb R^n)
$$
so $F\in \widetilde C^{\infty, c}_{N-\ell}(\mathbb R^N)$ and
$$\widetilde T_{\ell,\aleph}F
\ =\ 
[\widetilde T_{\ell,\aleph}F]_N\,dz_1\wedge\cdots\wedge dz_N\ \in\  \widetilde C^{\infty, c}_N(\mathbb R^N).$$
 Using \eqref{E:sobolev-deriv-3} and \eqref{E:def-F_I} we find
\begin{eqnarray*}
[\widetilde T_{\ell, \aleph} F]_N\circ i = 
\left(
\sum\limits_{\alpha\in \mathcal S(n, k)}
\epsilon^{\aleph(i\alpha)\,\aleph(i\alpha)'}_{1,\ldots, N}
\frac{\dee^k}{\dee z^{i\alpha}}F_{\aleph(i\alpha)'}
\right)\circ i
=\\
=\sum\limits_{\alpha\in \mathcal S(n, k)}
\frac{\dee^kg_\alpha}{\dee x^{\alpha}} = 0.
\end{eqnarray*}
where the last identity is due to the hypothesis 
\eqref{E:VS-div-k}.
Now observe that if $G\in \widetilde\Lambda_q(\mathbb R^N)$ then
$$
G=0 \iff G_I\circ i =0 \ \ \ \text{for\ each}\ \ I\in\mathcal I(N, q).
$$
 Combining all of the above we obtain
 $$
 \widetilde T_{\ell,\aleph}F=0
 $$
 so that Theorem \ref{P:LR-inflated-dim} grants
 \begin{equation*}
 |\langle F, H\rangle_{\widetilde L}|\leq \widetilde C
 \|F\|_{\widetilde L^1_q(\mathbb R^N)}
 \|\nabla H\|_{\widetilde L^n_q(\mathbb R^N)}.
 \end{equation*}
But since $(\pi\circ i) (x) =x$ for all $x\in \mathbb R^n$ it follows from 
 \eqref{E:pairing-big-1}, 
  \eqref{E:def-F_I} and
 \eqref{E:def-H_I} that
$$
\langle F, H\rangle_{\widetilde L} = \!\!\!\sum\limits_{\alpha\in\mathcal S(n, k)}\int
  g_\alpha\cdot h_\alpha;\  \|F\|_{\widetilde L^1_q(\mathbb R^N)} = 
  \|g\|_{L^1};\  \|\nabla H\|_{\widetilde L^n_q(\mathbb R^N)} =
  \|\nabla h\|_n
$$
and so
\begin{equation}\label{E:technical}
  \left|\,\sum\limits_{\alpha\in\mathcal S(n, k)}\int
  g_\alpha\cdot h_\alpha\ \right|\leq 
   \widetilde C\|g\|_{L^1}\,\|\nabla h\|_{L^n}
\end{equation}
is true for any $h\in (L^\infty\cap W^{1, n})(\mathbb R^n, \mathbb R^m)$. Now 
fix $\alpha_0\in \mathcal S(n, k)$ arbitrarily and define
$$
\hat h_\alpha := \delta_{\alpha_0\alpha}\, h_\alpha,\quad \alpha\in \mathcal S(n, k)
$$
where $\delta_{\alpha_0\alpha}$ denotes the Kroenecker symbol.
Then $\hat h\in (L^\infty\cap W^{1, n})(\mathbb R^n, \mathbb R^m)$ and so
 by applying \eqref{E:technical} to $\hat h$ we obtain
$$
\left|\sum\limits_{\alpha\in\mathcal S(n, k)}
\int\limits g_\alpha\cdot \hat h_\alpha \right|\ \leq\  \widetilde C\, \|g\|_{L^1}\
|\nabla \hat h\|_{L^n}.
$$
However
$$
\left|\sum\limits_{\alpha\in\mathcal S(n, k)}
\int\limits g_\alpha\cdot \hat h_\alpha \right| =
\left|
\int\limits g_{\alpha_0}\cdot h_{\alpha_0} \right|\quad \mbox{and}\quad
\|\nabla\hat h\|_{L^n}\leq \|\nabla h\|_{L^n},
$$
so \eqref{E:VS-ineq-div-k} is true for any choice of $\alpha_0\in\mathcal S(n, k)$, with $C:=\widetilde C$.
\qed

\bigskip

 {\it {Theorem \ref{P:LR-inflated-dim}\ $\Rightarrow$\ 
 Theorem \ref{P:LR-source-smallest}.}}\quad 
   Let $\ell$ be an admissible degree increment such that 
   $$
   n\ge \ell,
   $$
   let $\aleph$ be any one-to-one correspondence: $i\mathcal S (n, k)\to \mathcal I(N,\ell)$  and let 
          $0 \leq q\leq n-\ell$ be given. 
          Suppose that 
    $$f=\sum\limits_{I\in\mathcal I(n, q)}\!\!\! f_I\,dx^I\in L^1_q(\mathbb R^n)$$ 
       satisfies the hypotheses of Theorem \ref{P:LR-source-smallest};
         without loss of generality we may assume that
      $f \in C^{\infty, c}_q(\mathbb R^n)$.
      By the definition of $\Top$, see \eqref{E:Top-def}, we have
        $$
      \Top f =
      \sum\limits_{L\in\mathcal I(n, q+\ell)}
      \left(
  \sum\limits_{\atopp{I\in\mathcal I(n, q)}{\aleph(i\alpha)\in\mathcal I(n, \ell)}}
  \!\!\!\!\!\!
  \epsilon^{\aleph(i\alpha) I}_L\left(\frac{\dee^k (f_I\circ \pi)}{\dee z^{i\alpha}}
\right)\circ i\right)\! dx^L
      $$
and applying \eqref{E:sobolev-deriv-3} we obtain
      \begin{eqnarray}\label{E:Top-comp}
      \Top f =
      \sum\limits_{L\in\mathcal I(n, q+\ell)}
      \left(
  \sum\limits_{\atopp{I\in\mathcal I(n, q)}{\aleph(i\alpha)\in\mathcal I(n, \ell)}}
  \!\!\!\!\!\!
  \epsilon^{\aleph(i\alpha) I}_L\,\frac{\dee^kf_I}{\dee x^\alpha}
\right)\! dx^L \  
=0\notag
    \end{eqnarray}
 where the last identity is due to the hypothesis 
 \eqref{E:LR-div-invar-k-source-smallest}.  
      Fix $I_0\in\mathcal I(n, q)$ and
    choose (any) $L_0\in\mathcal I(n, q+\ell)$ so that
    \begin{equation*}
    I_0\subseteq L_0.
    \end{equation*}
    (The hypothesis $q\leq n-\ell$ 
       grants $q+\ell \leq  n$, so at least one such $L_0$ must exist.)
    
    Note that since $\ell\leq n\leq N$ we have $\mathcal I(n, \ell)\subseteq \mathcal I(N, \ell)$, so with $I_0$ and $L_0$ fixed as above, we may define 
    $$F^{L_0}=\sum\limits_{J\in\mathcal I(N, \ell)}\!\!\!
    F_J^{L_0}\, dz^J\in \widetilde C^{\infty, c}_\ell(\mathbb R^N)$$
     via
    \begin{equation*}
    F_J^{L_0}:=\displaystyle{\left\{\begin{array}{rcl}
    \sum\limits_{I\in\mathcal I(n, q)}\!\!\!\epsilon^{J I}_{L_0}\ f_I\circ\pi 
    &\mbox{for}& J\in \mathcal I(n, \ell)\\
    0&\mbox{for}& J\in\mathcal I(N, \ell)\setminus \mathcal I(n, \ell).
    \end{array}
    \right.}
    \end{equation*}
Applying \eqref{E:def-hilbert-tilde-T-adj} with $q:=\ell$ we obtain (ignore the factor $(-1)^{k+ q(N-\ell-q)}$) 
\begin{equation*}
\widetilde T^*_{\ell, \aleph} F^{L_0} = \sum\limits_{\alpha\in\mathcal S(n, k)}\!\!
\frac{\dee^k F^{L_0}_{\aleph (i\alpha)}}{\dee z^{i\alpha}}\in \widetilde C^{\infty, c}_0(\mathbb R^N)
\end{equation*}
and by the definition of $F^{L_0}$ this is further simplified to
$$
\widetilde T^*_{\ell, \aleph} F^{L_0} =
\sum\limits_{\atopp{I\in\mathcal I(n, q)}{\aleph(i\alpha)\in\mathcal I(n, \ell)}}
  \!\!\!\!\!\! \epsilon^{\aleph (i\alpha) I}_{L_0}\,\frac{\dee^k(f_I\circ \pi)}
  {\dee z^{i\alpha}}.
$$
 Note that on account of \eqref{E:pi-i} and of \eqref{E:sobolev-deriv-3} we have
$$
\Big(\widetilde T^*_{\ell, \aleph} F^{L_0}\Big)\circ i =
\sum\limits_{\atopp{I\in\mathcal I(n, q)}{\aleph(i\alpha)\in\mathcal I(n, \ell)}}
  \!\!\!\!\!\! \epsilon^{\aleph (i\alpha) I}_{L_0}\,\frac{\dee^k f_I}
  {\dee x^{\alpha}} = [\Top f]_{L_0} =0.
$$
But 
$\widetilde T^*_{\ell, \aleph} F^{L_0}\in \widetilde \Lambda_0(\mathbb R^N)$ and so
$$
\Big(\widetilde T^*_{\ell, \aleph} F^{L_0}\Big)\circ i = 0 \iff 
\widetilde T^*_{\ell, \aleph} F^{L_0} = 0,
$$
see \eqref{E:def-tilde-q-N}. Thus
$$
\widetilde T^*_{\ell, \aleph} F^{L_0} = 0
$$
and by Theorem \ref{P:LR-inflated-dim}  we conclude that
\begin{equation}\label{E:aux3}
|\langle F^{L_0}, H\rangle_{\widetilde L}|\leq \widetilde C
\| F^{L_0}\|_{\widetilde L^1_\ell}\|\nabla H\|_{\widetilde L^n_\ell}
\end{equation}
is true for any $H\in \widetilde C^{\infty, c}_\ell(\mathbb R^N)$. Now set
$$
J_0:= L_0\setminus I_0\in\mathcal I(n, \ell)
$$
and let  
 $$h=\sum\limits_{I\in\mathcal I(n, q)}\!\!\!h_I\,dx^I\in (L^\infty_q\cap 
     W^{1, n}_q)(\mathbb R^n)$$ 
be given (without loss of generality we may assume that $h\in C^{\infty, c}_q(\mathbb R^n)$).
Define $$\hat H = \sum\limits_{J\in\mathcal I(N, \ell)}\!\!\hat H_J\,dz^J\in
\Lambda_\ell(\mathbb R^N)$$ with
$$
\hat H_J =\delta_{J_0 J}\!\!\!\sum\limits_{I\in\mathcal I(n,  q)}\!\!\epsilon^{JI}_{L_0}\,h_I\circ \pi,\quad J\in\mathcal I(N,\ell),
$$
where $\delta_{J_0J}$ is the Kroenecker symbol. Then
$$
\hat H_J\circ i\circ \pi := \hat H_J\quad\mbox{and}\quad \hat H_J\circ i\in C^\infty_c(\mathbb R^n),\quad J\in\mathcal I(N, \ell)
$$
so that 
$$\hat H\in\widetilde C^{\infty, c}_\ell (\mathbb R^N).$$
Note that 
$$
\langle F^{L_0}, \hat H\rangle_{\widetilde L} = \int\limits_{\mathbb R^n}
\! f_{I_0}\cdot h_{I_0},\quad\mbox{and}\quad 
\|\nabla \hat H\|_{\widetilde L^n(\mathbb R^N)} \lesssim
 \|\nabla h\|_{L^n(\mathbb R^n)}.
$$
Moreover, by the definition of $F^{L_0}$ we have
$$
\|F^{L_0}\|_{\widetilde L^1_\ell(\mathbb R^N)}\lesssim
\|f\|_{L^1(\mathbb R^n)}.
$$
Thus, applying \eqref{E:aux3} to $\hat H$ we conclude that
\begin{equation}\label{E:aux4}
\left | \, \int\limits_{\mathbb R^n}\!\! f_{I_0}\cdot h_{I_0}\right| \leq
\widetilde C \|f\|_{L^1(\mathbb R^n)}\|\nabla h\|_{L^n(\mathbb R^n)}
\end{equation}
is true for any $I_0\in\mathcal I(n, q)$, for any $0\leq q\leq n-\ell$ and for any
$h\in C^{\infty, c}_q(\mathbb R^n)$, and this in turn implies \eqref{E:LR-ineq-source-smallest}.

In order to prove \eqref{E:LR-ineq-source-smallest-adj}, it suffices to apply
\eqref{E:LR-ineq-source-smallest} to: $f:= *_n g$ (with $q:=n-p$).
 \qed
  \bigskip

 {\it {Theorem \ref{P:LR-source-smallest}\ $\Rightarrow$\ 
 Theorem
  \ref{T:BB-thm-3}.}}\quad We claim that, in fact, Theorem \ref{T:BB-thm-3}
  is equivalent to the statement for $\AdjTop$ in Theorem \ref{P:LR-source-smallest} in the special case:
  $\ell =1$; $q=1$ and for specific choices of the ordering $\aleph$.
  Indeed it is easy to see that, for 
  $\ell =1$ and $q=1$, \eqref{E:Top-comp} and \eqref{E:def-Top-adj} give
  $$
  \mathcal T_{1, \aleph}^*f=
  (-1)^{k+n}\sum\limits_{j=1}^n\frac{\dee^k f_j}{\dee x^{\pi\circ \aleph^{-1}(j)}},
  \quad
 f=\sum\limits_{j=1}^n
 f_{j}\,dx_{j}\in \Lambda_{1}(\mathbb R^n).
  $$
  Choosing now any ordering $\aleph: i\mathcal S(n, k)\leftrightarrow \mathcal I(N, 1)$ 
  such that 
  $$
  \pi\circ\aleph^{-1} (j) = (0, \ldots, 0, k,0,\ldots, 0), \quad j=1,\ldots, n
  $$
 (where, in the expression above, it is understood that $k$ occupies the $j$-th position) we obtain
  $$
  \mathcal T_{1, \aleph}^*f =
  \sum\limits_{j=1}^n\frac{\dee^ku_j}{\dee x_j^k}, \quad u_j := (-1)^{k+n}f_j , \ j=1, \ldots, n.
  $$
  The equivalence of the two statements is now apparent.
\qed
\bigskip

{\it Proof of Lemma \ref{L:box-reps}}.\quad The proof of 
\eqref{E:big-box-ell-repr} and \eqref{E:big-coeffs} is a computation that uses \eqref{E:gen-LR-oper} along with the following  coordinate-based representation for $\widetilde T^*_{\ell, \aleph}$, which is obtained from 
\eqref{E:def-hilbert-tilde-T-adj}:
$$
\widetilde T^*_{\ell, \aleph}H = (-1)^{k+N\ell}\!\!\!\!\!
\sum\limits_{\atopp{V\in\mathcal I(N, q-\ell)}
{\atopp{I\in\mathcal I(N, q)}{\beta\in\mathcal S(n, k)}}}
\!\!\!\!\!\epsilon_I^{\aleph(i\beta)V}\,
\frac{\dee^kH_I}{\dee z^{i\beta}}\,dz^V
$$
To prove \eqref{E:big-box-1-repr} we examine
\eqref{E:big-coeffs} in the case $\ell =1$:
\begin{equation}
\widetilde C^{MI}_{\aleph(i\alpha)\,\aleph(i\beta)} =
 \widetilde A^{MI}_{\aleph(i\alpha)\,\aleph(i\beta)}+ 
 \widetilde B^{MI}_{\aleph(i\alpha)\,\aleph(i\beta)}
\end{equation}
where
\begin{equation}\label{E:def-A}
 \widetilde A^{MI}_{\aleph(i\alpha)\,\aleph(i\beta)}\ := 
 \sum\limits_{L\in\mathcal I(N, q+1)}
\!\!\!\!\!\! \epsilon^L_{\aleph(i\alpha) I}\cdot\epsilon^L_{\aleph(i\beta)M}\, ,
\end{equation}
\begin{equation}\label{E:def-B}
\widetilde B^{MI}_{\aleph(i\alpha)\,\aleph(i\beta)}\ :=
\!\!\!\sum\limits_{K\in\mathcal I(N, q-1)}
\!\!\!\!\!\! \epsilon^M_{\aleph(i\alpha) K}\!\cdot\epsilon^I_{\aleph(i\beta) K}
\end{equation}
and distinguish two cases: $\alpha \neq \beta$; and $\alpha = \beta$.
\medskip

 \noindent Suppose first that $\alpha\neq \beta$. In this case we claim that 
 $\widetilde C^{MI}_{\aleph(i\alpha)\, \aleph(i\beta)}=0$.
 The proof of this claim rests on the following 
 \begin{rem}\label{R:rem}
 The truth 
 value of the following three (combined) conditions on $\aleph,\ \alpha,\ \beta,\ I$ and $M$:
 \begin{equation}\label{E:condition-1}
 \aleph(i\alpha)\notin \{I\};\  \aleph(i\beta)\notin \{M\};\ \
 \{\aleph(i\alpha)\}\cup \{I\} = \{\aleph(i\beta)\}\cup \{M\}
 \end{equation}
 
 is equivalent\footnote{If $q = 0$ or $q=N-1$ then \eqref{E:condition-1} is equivalent to  \eqref{E:condition-2} in the sense that each is false.}
 to the truth value of \begin{equation}\label{E:condition-2}
 \aleph(i\alpha)\in \{M\};\ \aleph(i\beta)\in\{I\}; \ \ \{M\}\setminus\{\aleph(i\alpha)\}=
 \{I\}\setminus\{\aleph(i\beta)\}.
 \end{equation}
 \end{rem}
 We postpone the proof of Remark \ref{R:rem} and continue with the proof of Lemma \ref{L:box-reps}; 
 to this end we claim  that if $\alpha\neq \beta$ then 
   $$
 \mbox{\eqref{E:condition-1}\ holds}\quad \iff\quad \widetilde A^{MI}_{\aleph(i\alpha)\aleph(i\beta)}
 \neq 0
 $$
 Indeed, since $\alpha,\ \beta, \ I$ and $M$ are fixed,  the summation
 that defines $\widetilde A^{MI}_{\aleph(i\alpha)\aleph(i\beta)} $, see\eqref{E:def-A}, consists of at most one term, that is
$$
\widetilde A^{MI}_{\aleph(i\alpha)\aleph(i\beta)} = \
\epsilon^{L_0}_{\aleph(i\beta)M}\!\cdot\epsilon^{L_0}_{\aleph(i\alpha) I} 
$$
for at most one choice of $L_0\in\mathcal I(N, q+1)$, and it's easy to see that
\eqref{E:condition-1} holds if, and only if, there is exactly one choice of  $L_0\in \mathcal I(N, q+1)$ such that $\epsilon^{L_0}_{\aleph(i\beta)M}\!\cdot\epsilon^{L_0}_{\aleph(i\alpha) I}\neq 0$ and in such 
case we have
  \begin{equation}\label{E:A-non-zero}
  \widetilde A^{MI}_{\aleph(i\alpha)\aleph(i\beta)} = \
  \ \epsilon^{\aleph(i\beta)M}_{\aleph(i\alpha) I}.
  \end{equation}
Similar considerations grant, again for $\alpha\neq\beta$, that we also have
$$
 \mbox{\eqref{E:condition-2}\ holds}\quad \iff\quad 
 \widetilde B^{MI}_{\aleph(i\alpha)\aleph(i\beta)}\neq 0
$$ 
and if $\widetilde B^{MI}_{\aleph(i\alpha)\aleph(i\beta)}\neq 0$ there is a unique choice of $K_0\in\mathcal I(N, q-1)$ such that 
\begin{equation}\label{E:B-non-zero}
  \widetilde  B^{MI}_{\aleph(i\alpha)\,\aleph(i\beta)} =
   \epsilon^M_{\aleph(i\alpha) K_0}\,\epsilon^I_{\aleph(i\beta) K_0}  =
 - \,\epsilon^{\aleph(i\beta)M}_{\aleph(i\alpha)I}\, .
    \end{equation}

Combining all of the above, we conclude that if $\alpha\neq\beta$ then either
$$
\widetilde A^{MI}_{\aleph(i\alpha)\aleph(i\beta)}\  = \ 
  \widetilde B^{MI}_{\aleph(i\alpha)\aleph(i\beta)}\ =\ 0
$$
or 
$$
\widetilde A^{MI}_{\aleph(i\alpha)\aleph(i\beta)}\  = \ 
 -\  \widetilde B^{MI}_{\aleph(i\alpha)\aleph(i\beta)}.
$$
In either case it follows that
\begin{equation}\label{E:big-C-1}
\widetilde C^{MI}_{\aleph(i\alpha)\aleph(i\beta)}\  = \ 0\quad 
\mbox{whenever}\ \alpha\neq \beta\, .
\end{equation}

Suppose next that  $\alpha = \beta$; in this case \eqref{E:def-A} and 
\eqref{E:def-B} become

\begin{equation}\label{E:def-A-1}
 \widetilde A^{MI}_{\aleph(i\alpha)\,\aleph(i\alpha)}\ = 
 \sum\limits_{L\in\mathcal I(N, q+1)}
\!\!\!\!\!\! \epsilon^L_{\aleph(i\alpha) I}\cdot\epsilon^L_{\aleph(i\alpha)M}\, ,
\end{equation}
\begin{equation}\label{E:def-B-1}
\widetilde B^{MI}_{\aleph(i\alpha)\,\aleph(i\alpha)}\ =
\!\!\!\sum\limits_{K\in\mathcal I(N, q-1)}
\!\!\!\!\!\! \epsilon^M_{\aleph(i\alpha) K}\!\cdot\epsilon^I_{\aleph(i\alpha) K}
\end{equation}
and since $\alpha,\ I$ and $M$ are fixed, each of these summations consists of
at most one term, that is
$$
\widetilde A^{MI}_{\aleph(i\alpha)\,\aleph(i\alpha)}\ = 
\epsilon^{L_0}_{\aleph(i\alpha) I}\cdot\epsilon^{L_0}_{\aleph(i\alpha)M}\, ;\quad
\widetilde B^{MI}_{\aleph(i\alpha)\,\aleph(i\alpha)}\ =
\epsilon^M_{\aleph(i\alpha) K_0}\!\cdot\epsilon^I_{\aleph(i\alpha) K_0}
$$
for at most one choice for each of $L_0\in\mathcal I(N, q+1)$ and
$K_0\in\mathcal I(N, q-1)$. In particular we see that
\begin{equation}
I\neq M\quad \Rightarrow\quad  \widetilde A^{MI}_{\aleph(i\alpha)\,\aleph(i\alpha)}\ =\ 
\widetilde B^{MI}_{\aleph(i\alpha)\,\aleph(i\alpha)}\ =\  0\, .
\end{equation}
On the other hand, for $I=M$ we have
$$
\widetilde A^{MM}_{\aleph(i\alpha)\,\aleph(i\alpha)}\ = 
(\epsilon^{L_0}_{\aleph(i\alpha) M})^2
\, ;\quad
\widetilde B^{MM}_{\aleph(i\alpha)\,\aleph(i\alpha)}\ =
(\epsilon^M_{\aleph(i\alpha) K_0})^2
$$
for at most one choice of $L_0$ and of $K_0$.
 We now further distinguish between $\aleph(i\alpha)\in \{I\}$
 and $\aleph(i\alpha)\notin \{I\}$. If $\aleph(i\alpha)\in \{I\}$ then we have
 $\widetilde A^{MM}_{\aleph(i\alpha)\,\aleph(i\alpha)} = 0$ (because the $L$'s do not have repeated terms) and $\widetilde B^{MM}_{\aleph(i\alpha)\,\aleph(i\alpha)} = 1$. If, instead, $\aleph(i\alpha)\notin \{I\}$ then we find by the same token that
 $\widetilde A^{MM}_{\aleph(i\alpha)\,\aleph(i\alpha)} = 1$ and 
 $\widetilde B^{MM}_{\aleph(i\alpha)\,\aleph(i\alpha)} = 0$. 
 All together this gives
 \begin{equation}\label{E:big-C-2}
 \widetilde C^{MI}_{\aleph(i\alpha)\,\aleph(i\alpha)} =
 \left\{ \begin{array}{rcl}
 0&\mbox{for}& M\neq I\\
 1 &\mbox{for}& M=I\, .
 \end{array}\right. 
 \end{equation}
Combining \eqref{E:big-C-1} and  \eqref{E:big-C-2} we obtain \eqref{E:big-box-1-repr}.
  The proofs of \eqref{E:box-ell-repr} -- \eqref{E:box-1-repr} are obtained in a similar fashion; in this case \eqref{E:def-hilbert-T-adj} grants
 $$
 \Top^*h = (-1)^{k+n\ell}\!\!\!\!\!\!\!\!\!\!
\sum\limits_{\atopp{V\in\mathcal I(n, q-\ell)}
{\atopp{I\in\mathcal I(n, q)}{\beta\in(\pi\circ\aleph^{-1})(\mathcal I(n, \ell))}}}
\!\!\!\!\!\!\!\!\!\!\epsilon_I^{\aleph(i\beta)V}\,
\frac{\dee^kh_I}{\dee x^\beta}\,dx^V.
$$
\qed

\bigskip

{\em Proof of Remark \ref{R:rem}.}\quad 
If $\alpha\neq \beta$ and the three conditions in \eqref{E:condition-1} hold, then it follows at once that the first two conditions in \eqref{E:condition-2}
are true;
 by the first condition in \eqref{E:condition-1} we have $\{I\} = \{I\}\cap
  \{\aleph(i\alpha)\}^c$; combining this identity with the third condition in \eqref{E:condition-1} we obtain
 $\{I\} = \left(\{\aleph(i\beta)\}\cup \{M\}\right)\cap \{\aleph(i\alpha)\}^c$, and since $\alpha\neq \beta$ then $\{\aleph (i\beta)\}\cap \{\aleph(i\alpha)\}^c =
 \{\aleph (i\beta)\}$, and it follows that 
 $$
 \{I\} =  \{\aleph (i\beta)\} \cup \{Q_0\},\quad \{Q_0\}:= \{M\}\cap \{\aleph(i\alpha)\}^c
 $$
 By the second condition in \eqref{E:condition-1} we have $\{\aleph (i\beta)\} \cap \{Q_0\}=\emptyset$ and so
 $$
 \{I\} \setminus \{\aleph (i\beta)\} = \{Q_0\}
 $$
 On the other hand, since we have proved that $\aleph (i\alpha)\in\{M\}$ is true, then we have 
 $$ 
 \{M\} = (\{M\}\cap \{\aleph (i\alpha)\}) \cup (\{M\}\cap \{\aleph (i\alpha)\}^c) =
 \{\aleph (i\alpha)\}\cup \{Q_0\}
 $$
 and obviously $\{\aleph (i\alpha)\}\cap \{Q_0\} =\emptyset$, so
 $$
 \{M\} \setminus  \{\aleph (i\alpha)\} = \{Q_0\}
 $$
 which shows that the third condition in \eqref{E:condition-2} is true, as well.
 
  Suppose, conversely, that $\alpha\neq \beta$ and that the three conditions in \eqref{E:condition-2} hold.
 Then the first condition in \eqref{E:condition-2} grants
 $\{M\} = \{\aleph (i\alpha)\}\dot{\cup} \{P_0\}$
 (where $\dot{\cup}$ denotes disjoint union); similarly, the second condition
 in \eqref{E:condition-2} grants $\{I\} = \{\aleph (i\beta)\}\dot{\cup} \{S_0\}$, and 
 it follows from the third condition in \eqref{E:condition-2} that $S_0=P_0$.
 Note that in particular $\aleph (i\alpha)\notin \{P_0\}$ and $\aleph (i\beta)\notin \{P_0\}$; since $\alpha\neq \beta$, it follows that the first two conditions in \eqref{E:condition-1} hold. But these (and the above) considerations in turn imply
 $$
 \{\aleph (i\alpha)\} \cup\{I\} = 
 \{\aleph (i\alpha)\}\cup  \{\aleph (i\beta)\} \cup P_0 =
  \{\aleph (i\beta)\} \cup\{M\}
 $$
 which shows that the third condition in \eqref{E:condition-1} is true, as well.
 \qed
\bigskip

{\em Proof of Lemma \ref{L:big-box-elliptic-estimates}.}\quad 
The proof is easily reduced to the classical theory via Corollary
\ref{C:Legendre-Hadamard} along with
\eqref{E:sobolev-deriv-3} and the coordinate-based representations for $\|\cdot \|_{\widetilde L^n}$, see 
\eqref{E:r-big-norm-1}. See \cite{CZ}, \cite[pg. 62]{SR}, \cite[VI.5]{S} and \cite[13.6]{T}.
\qed
\bigskip

{\em Theorem \ref{P:LR-inflated-dim} ($\ell =1$) $\Rightarrow$\ 
Theorem \ref{T:LR-thm-Z}.}\quad Without loss of generality we may assume: 
$F\in\widetilde C^{\infty, c}_{q+1}(\mathbb R^N)$; 
$G\in\widetilde C^{\infty, c}_{q-1}(\mathbb R^N)$, so that
$Z\in\widetilde C^{\infty, c}_{q}(\mathbb R^N)$.
Write
$$
Z=X + Y
$$
where
\begin{equation}\label{E:system-F}
\left\{
\begin{array}{clr}
\widetilde T_{1, \aleph} X&=&F\\
\widetilde T_{1, \aleph}^* X&=&0
\end{array}
\right.
\end{equation}
and 
\begin{equation}\label{E:system-G}
\left\{
\begin{array}{clr}
\widetilde T_{1, \aleph} Y&=&0\\
\widetilde T_{1, \aleph}^* Y&=&G
\end{array}
\right.
\end{equation}
We claim that
\begin{equation}\label{E:F-est}
\|X\|_{\widetilde W^{k-1, r}_q}\leq C \|F\|_{\widetilde L^1_{q+1}},\quad r:=n/(n-1)
\end{equation}
and 
\begin{equation}\label{E:G-est}
\|Y\|_{\widetilde W^{k-1, r}_q}\leq C \|G\|_{\widetilde L^1_{q-1}},\quad r:=n/(n-1)
\end{equation}
Note that if $Y$ solves \eqref{E:system-G} then $X:= *_NY$ solves 
\eqref{E:system-F} with $F:= *_NG$, and so it suffices to prove 
\eqref{E:F-est} for $F$ and $X$ as in \eqref{E:system-F}. 
By duality, this is equivalent to proving
\begin{equation}\label{E:F-est-aux}
\big|\langle D^{i\beta} X, \varphi\rangle_{\widetilde L}\big|\leq
C\|F\|_{\widetilde L^1_{q+1}} \|\varphi\|_{\widetilde L_q^{n}}
\end{equation}
for any $\beta\in\mathcal S(n, s)$ with $0\leq s\leq k-1$, and for any $\varphi\in \widetilde C^{\infty, c}_q(\mathbb R^N)$. Let $\Phi\in \widetilde C^{\infty, c}_q(\mathbb R^N)$ be as in Lemma \ref{L:big-box-elliptic-estimates}. Note that
$$
\widetilde T_{1, \aleph}D^{i\beta}X = D^{i\beta}\,\widetilde T_{1, \aleph}X =
 D^{i\beta} F
;\qquad 
\widetilde T^*_{1, \aleph}D^{i\beta}X = D^{i\beta}\,\widetilde T^*_{1, \aleph}X = 0.
$$
By \eqref{E:hilbert-tilde-T-adj} and the above considerations it follows that
$$
\langle D^{i\beta}X, \varphi\rangle_{\widetilde L} =
\langle D^{i\beta}F, \widetilde T_{1, \aleph}\Phi\,\rangle_{\widetilde L}  =
\langle F, D^{i\beta}\widetilde T_{1, \aleph}\Phi\,\rangle_{\widetilde L}  
$$
By Theorem \ref{P:LR-inflated-dim}  ($\ell =1$; $H:= D^{i\beta}\widetilde T_{1, \aleph}\Phi \in \widetilde C^{\infty, c}_{q+1}(\mathbb R^N)$) we have
$$
\big| \langle D^{i\beta}X, \varphi\rangle_{\widetilde L}\big| \leq 
C\|F\|_{\widetilde L^1_{q+1}}\|\nabla
 D^{i\beta}\widetilde T_{1, \aleph}\Phi\|_{\widetilde L^n_{q+1}} \leq
 C\|F\|_{\widetilde L^1_{q+1}}\|\Phi\|_{\widetilde W_{q}^{2k, n}},
$$
and it follows from Lemma \ref{L:big-box-elliptic-estimates} (with $p:= n$) that
$$
\big| \langle D^{i\beta}X, \varphi\rangle_{\widetilde L}\big| \leq
 C\|F\|_{\widetilde L^1_{q+1}}\|\varphi\|_{\widetilde L_{q}^{ n}}
$$
as desired.
\qed
\bigskip

{\em Theorem \ref{T:LR-thm-Z}\ $\Rightarrow$\ 
Theorem \ref{P:LR-inflated-dim} ($\ell =1$;\quad $1\leq q, \, p\leq N-1$).}\quad
\smallskip

\noindent We show that \eqref{E:LR-ineq-inflated} holds for any $q$ in the range
$1\leq q\leq N-1$. Suppose that $\widetilde T_{1, \aleph} F=0$, 
$F\in \widetilde L^1_{q}(\mathbb R^N)$ and let 
$H\in (\widetilde L^\infty_{q}\cap \widetilde W^{1, n}_q)(\mathbb R^N)$.
Without loss of generality we may assume: 
$H, F\in \widetilde C^{\infty, c}_{q}(\mathbb R^N)$. Let 
$X\in\widetilde C^{\infty, c}_{q-1}(\mathbb R^N)$ be the solution of \eqref{E:system-F} with data $F$.
Then, by H\" older inequality \eqref{E:holder-big} we have
$$
\big|\langle F, H\rangle_{\widetilde L}\big| = 
\big |\langle X, \widetilde T_{1, \aleph}^* H\rangle_{\widetilde L}\big |
\lesssim
\|X\|_{\widetilde W^{k-1, r}_q}
\|\widetilde T_{1, \aleph}^* H\|_{\widetilde W^{-(k-1), n}_q}
$$
and it follows from the latter and Theorem  \ref{T:LR-thm-Z} that
$$
\big|\langle F, H\rangle_{\widetilde L}\big|
\lesssim \|F\|_{\widetilde L^1_{q}}
\|\widetilde T_{1, \aleph}^* H\|_{\widetilde W^{-(k-1), n}_{q-1}}
$$
Now observe that if integrate the expression 
$$\langle \widetilde T_{1, \aleph}^* H, \zeta \rangle_{\widetilde L}$$ 
by parts $(k-1)$-times and then apply 
 H\"older inequality we obtain
 $$
|\langle \widetilde T_{1, \aleph}^* H, \zeta \rangle_{\widetilde L}| 
\leq \|\nabla H\|_{\widetilde L^n_{q}}\|\zeta\|_{\widetilde W^{k-1, r}_{q-1}}
$$
and this leads to the conclusion of the proof of \eqref{E:LR-ineq-inflated} as 
$$
\|\widetilde T_{1, \aleph}^* H\|_{\widetilde W^{-(k-1), n} _{q-1}} =
\sup\limits_{\|\zeta\|_{\widetilde W^{k-1, r}_{q-1}}\,\leq 1}
\big |\langle \widetilde T_{1, \aleph}^*H, \zeta
 \rangle_{\widetilde L}\big |
$$
\qed
\bigskip

\section{Concluding Remarks}\label{S:open-problems}
\medskip

  \begin{itemize}
  \item[1.] The proof of Theorems \ref{P:LR-inflated-dim} and 
  \ref{T:LR-thm-Z} rely on the specific choice of the embedding
  $i: \mathbb R^n\to \mathbb R^N$ only to the extent that \eqref{E:pi-i},
  in fact \eqref{E:pi-i-pr}, and \eqref{E:sobolev-deriv-2} hold. This suggests that 
  Theorems \ref{P:LR-inflated-dim}, \ref{P:LR-source-smallest} and 
  \ref{T:LR-thm-Z} should hold in the more general context of an isometrically embedded manifold
 $$\mathcal M^{(n)}\hookrightarrow \mathbb R^N$$

  \item[]
  \item[2.] If $q\ge N-\ell+1$ or $p\leq \ell-1$ then one of the two
 compatibility conditions \eqref{E:LR-div-invar-k-inflated} 
 and \eqref{E:tilde-T-star-LR} holds trivially
  and in this case the conclusion of Theorems \ref{P:LR-inflated-dim} and \ref{T:LR-thm-Z} are  easily seen to be 
 false: if $k=1$ and $\widetilde T_{1,\aleph_0}=d$ (exterior derivative) substitute inequalities hold provided the 
  ``defective''
 data belongs to a suitable (proper) subspace of $L^1$, namely the real Hardy space 
 $H^1 (\mathbb R^n)$, see \cite{LS}. We do not know whether substitute inequalities 
   hold when $k\ge 2$.
 \item[]
 \item[3.]   In the context of \cite{VS4} our results say the following:
\item[]
\medskip

\begin{itemize}
 \item[$\centerdot$] $\widetilde T_{1, \aleph}$ is 
  canceling from
 $V:= \widetilde C^{\infty, c}_q(\mathbb R^N)$ to 
 $E:=  
 \widetilde C^{\infty, c}_{q+1}(\mathbb R^N)$ for any $0\leq q\leq N-2$, see \cite[Theorem 1.3]{VS4}.
 \item[] 
 
  \item[$\centerdot$] $\widetilde T^*_{1, \aleph}$ is 
   canceling from
 $V:= \widetilde C^{\infty, c}_q(\mathbb R^N)$ to
  $E:=  \widetilde C^{\infty, c}_{q-1}(\mathbb R^N)$ for any $2\leq q\leq N$, see \cite[Theorem 1.3]{VS4}.
  \item[] 
  
    \item[$\centerdot$] For any admissible degree increment $\ell$ and for any $0\leq q\leq N-\ell$, 
  $\widetilde T_{\ell, \aleph}$ is cocanceling from 
  $V:= \widetilde C^{\infty, c}_q(\mathbb R^N)$ to $E:= \widetilde C^{\infty, c}_{q+\ell}(\mathbb R^N)$, see
  \cite[Propositions 2.1 and 2.2]{VS4}.
  \item[] 
  
  \item[$\centerdot$] For any admissible degree increment $\ell$ and  for any 
  $\ell\leq q\leq N$, 
  $\widetilde T^*_{\ell, \aleph}$ is cocanceling from 
  $V:= \widetilde C^{\infty, c}_q(\mathbb R^N)$ to 
  $E:= \widetilde C^{\infty, c}_{q-\ell}(\mathbb R^N)$, see
   \cite[Propositions 2.1 and 2.2]{VS4}.
 \item[]
 
\item[$\centerdot$]  The class $\mathcal T_{\ell, \aleph}$ has similar properties with $V= C^\infty_q(\mathbb R^n)$ and $E = C^{\infty}_{q\pm\ell}(\mathbb R^n)$.
\end{itemize}
\medskip

\noindent In particular, $\widetilde T_{1, \aleph}$ and $\mathcal T_{1, \aleph}$ as well as their adjoints,  are new examples
of canceling operators of {\em arbitrary} order $k$.
 \end{itemize}

\end{document}